\documentclass[12pt,twoside,letter]{article}
\usepackage{amsmath,amssymb,euscript}

\title{Effective Dynamics of Magnetic Vortices}
\author{S. Gustafson, \quad I.M. Sigal}

\date{\today}

\newtheorem{thm}{Theorem}
\newtheorem{lem}{Lemma}
\newtheorem{prop}{Proposition}
\newtheorem{rem}{Remark}

\newcommand{\donothing}[1]{}

\begin{document}

%%%%%%%%%%%%%%%%%%%%%%%%%%%%%%%%%%%%%%%
%%%%%%%       shortcuts         %%%%%%%
%%%%%%%%%%%%%%%%%%%%%%%%%%%%%%%%%%%%%%%

\newcommand{\nc}{\newcommand}

\nc{\be}{\begin{equation}}
\nc{\la}{\label}
\nc{\ba}{\begin{array}}
\nc{\ea}{\end{array}}
\nc{\bs}{\begin{split}}
\nc{\es}{\end{split}}

\nc{\R}{\mathbb R}
\nc{\Z}{\mathbb Z}
\nc{\C}{\mathbb C}
\nc{\J}{\mathbb J}

\nc{\p}{\partial}
\nc{\pt}{\partial_t}
\nc{\ptt}{\partial_t^2}

\nc{\zb}{\underbar{z}}
\nc{\pb}{\underbar{p}}
\nc{\bb}{\underbar{b}}
\nc{\zbt}{\underbar{z}(t)}

\nc{\dA}{\nabla_A}

\nc{\nj}{(n_j)}
\nc{\nk}{(n_k)}
\nc{\nl}{(n_l)}
\nc{\nr}{(n_r)}
\nc{\nt}{(n_t)}

\nc{\vz}{v_{\zb,\chi}}
\nc{\vzo}{v_{\zb_0,\chi_0}}
\nc{\vzt}{v_{\zb(t),\chi(t)}}
\nc{\psiz}{\psi_{\zb,\chi}}
\nc{\Az}{A_{\zb,\chi}}
\nc{\Bz}{B_{\zb,\chi}}
\nc{\jz}{j_{\zb,\chi}}
\nc{\Tz}{T^{(\zb,\chi)}}
\nc{\Gz}{G^{(\zb,\chi)}}
\nc{\Lz}{L_{\zb,\chi}}
\nc{\Lzo}{L_{\zb_0,\chi_0}}
\nc{\Pz}{P_{\zb,\chi}}
\nc{\bPz}{\bar{P}_{\zb,\chi}}
\nc{\phiz}{\phi_{\zb,\chi,\pb,\zeta}}
\nc{\phizo}{\phi_{\zb_0,\chi_0,\pb_0,\zeta_0}}
\nc{\phizt}{\phi_{\zb(t),\chi(t),\pb(t),\zeta(t)}}
\nc{\wz}{w_{\zb,\chi,\pb,\zeta}}
\nc{\wzo}{w_{\zb_0,\chi_0,\pb_0,\zeta_0}}

\nc{\E}{{\cal E}}
\nc{\cH}{{\cal H}}

\nc{\e}{\epsilon}
\nc{\lam}{\lambda}
\nc{\G}{\Gamma}
\nc{\g}{\gamma}
\nc{\al}{\alpha}
\nc{\del}{\delta}
\nc{\Om}{\Omega}
\nc{\Omt}{\tilde{\Omega}}
\nc{\ta}{\tau}
\nc{\w}{\omega}
\nc{\io}{\iota}
\nc{\h}{\theta}
\nc{\z}{\zeta}
\nc{\s}{\sigma}
\nc{\Si}{\Sigma}
\nc{\Lam}{\Lambda}

\nc{\bP}{\bar{P}}
\nc{\bQ}{\bar{Q}}
\nc{\bL}{\bar{L}}

\nc{\chit}{\tilde{\chi}}
\nc{\It}{\tilde{I}}
\nc{\alt}{\tilde{\alpha}}
\nc{\Mt}{\tilde{M}}
\nc{\gt}{\tilde{\gamma}}
\nc{\Jt}{\tilde{J}}

\nc{\ra}{\rightarrow}

\nc{\ran}{\rangle}
\nc{\lan}{\langle}

\nc{\bfone}{{\bf 1}}

%%%%%%%%%%%%%%%%%%%%%%%%%%%%%%%%%%%%%%%%%%%
%%%%%%     Title and abstract      %%%%%%%%
%%%%%%%%%%%%%%%%%%%%%%%%%%%%%%%%%%%%%%%%%%%

\maketitle

\begin{abstract}
We study solutions of Ginzburg-Landau-type
evolution equations (both dissipative and Hamiltonian)
with initial data representing collections of
widely-spaced vortices. We show that for long times,
the solutions continue to describe collections of
vortices, and we identify (to leading order in the
vortex separation) the dynamical system describing
the motion of the vortex centres ({\it effective dynamics}).
%Since the basic techniques are largely equation-independent,
%we also outline a general theory of {\it effective dynamics}
%of ``approximate static solutions'' of
%an abstract evolution equation.  The results for Ginzburg-Landau
%equations are applications of these general ideas.
\end{abstract}

\tableofcontents

%%%%%%%%%%%%%%%%%%%%%%%%%%%%%%%%%%%%%%%%%%
%%%%%%     Introduction           %%%%%%%%
%%%%%%%%%%%%%%%%%%%%%%%%%%%%%%%%%%%%%%%%%%

\section{Introduction}

In this paper we study effective dynamics of magnetic (Abrikosov)
vortices in a macroscopic model of superconductivity, and of
Nielsen-Olesen or Nambu strings in the Abelian Higgs model of
particle physics. In both cases the equilibrium configurations are
described by the Ginzburg-Landau equations:
\begin{equation}
\la{eq:gl}
   \ba{c}
   -\Delta_A \psi = \lambda (1-|\psi|^2) \psi   \\
   curl^2 A = Im ( \bar{\psi} \nabla_A \psi )
   \ea
\end{equation}
where $(\psi , A) : \R^2 \to \C \times \R^2$,
$\nabla_A = \nabla - iA$, and $\Delta_A = \nabla^2_A$,
the covariant derivative and covariant Laplacian,
respectively.  Equations~(\ref{eq:gl})
are the Euler-Lagrange equations for the
Ginzburg-Landau energy functional
\be
\la{eq:en}
  \E_{GL}(\psi,A) := \frac{1}{2} \int_{\R^2} \left\{
  |\nabla_A \psi|^2 + (curl A)^2
  + \frac{\lam}{2}(|\psi|^2-1)^2 \right\}.
\end{equation}
In the case of superconductivity, the function $\psi : \R^2 \to
\C$ is called the {\it order parameter}; $|\psi |^2$ gives the density
of superconducting electrons.  The vector field $A : \R^2 \to
\R^2$ is the magnetic potential. The r.h.s. of the equation for
$A$ is the superconducting current. In the case of particle
physics, $\psi$ and $A$ are the Higgs and Abelian gauge
(electro-magnetic) fields, respectively.
(See \cite{rit} for reviews, and \cite{no}
for historical and physics background.)

In addition to being translationally
and rotationally invariant,
equations~(\ref{eq:gl}) are invariant under
gauge transformations:
\[
  (\psi,A) \mapsto (e^{i\chi}\psi, A + \nabla \chi)
\]
for any $\chi : \R^2 \ra \R$
(solutions are mapped to solutions under this
transformation).

We consider various time-dependent
versions of the Ginzburg-Landau equations~(\ref{eq:gl}).
The first example is the gradient-flow equations
\begin{equation}
\la{eq:gf}
\bs
   \partial_t \psi &= \Delta_A \psi
   + \lambda (1-|\psi|^2) \psi \\
   \partial_t A &= -curl^2 A + Im ( \bar{\psi} \nabla_A \psi ),
\end{split}
\end{equation}
a model in superconductivity theory (\cite{ge,t}).
We will refer to equations~(\ref{eq:gf}) as the
{\it superconductor model}
(they are sometimes called the
{\it Gorkov-Eliashberg equations} or
{\it time-dependent Ginzburg-Landau equations}).

The second example is
\begin{equation}
\la{eq:mh}
\bs
   \partial_t^2 \psi &= \Delta_A \psi
   + \lambda (1-|\psi|^2) \psi \\
   \partial_t^2 A &= -curl^2 A + Im ( \bar{\psi} \nabla_A \psi ),
\end{split}
\end{equation}
coupled (covariant) wave equations
describing the $U(1)$-gauge Higgs model of
elementary particle physics (\cite{jt})
(written here in the {\it temporal gauge}).
We will refer to equations~(\ref{eq:mh}) as
the {\it Higgs model} (they are sometimes also
called the {\it Maxwell-Higgs equations}).

The general framework we develop in this paper
also applies to coupled (complex)
Schr\"odinger and Maxwell equations
\begin{equation}
\label{eq:schmax}
\begin{split}
  & \g \p_t \psi = \Delta_A \psi + \lam(1-|\psi|^2)\psi  \\
  & \p_t^2 A = -curl^2 A + Im(\bar{\psi} \nabla_A \psi)
\end{split}
\end{equation}
with $Re \g \geq 0$,
or the Chern-Simons variant of these equations,
though the implementation for $Re \g = 0$
requires some additional technical steps.

Finite energy states $(\psi,A)$ are classified by the
topological degree
\[
  deg(\psi) := deg
  \left( \left. \frac{\psi}{|\psi|} \right|_{|x|=R} \right),
\]
where $R$ is sufficiently large (the winding
number of $\psi$ at infinity).  For each such
state we have the quantization of magnetic flux:
\[
  \int_{\R^2} B = 2 \pi deg(\psi) \in 2 \pi \Z,
\]
where $B := curl A$ is the magnetic field associated
with the vector potential $A$.

In each case the equations have ``radially symmetric''
(more precisely {\it equivariant}) solutions of the form
\begin{equation}
\label{eq:vort}
   \psi^{(n)} (x) = f_n (r) e^{in\theta} {\hbox{\quad and \quad}}
   A^{(n)}(x) = a_n (r) \nabla (n\theta) \ ,
\end{equation}
where $n$ is an integer and $(r,\theta)$ are the polar
coordinates of $x \in \R^2$.  As $r \to \infty$,
$a_n(r)$ and $f_n(r)$ converge to $1$ exponentially
fast with the rates $1$ and
$m_\lam := \min(\sqrt{2\lambda},2)$, respectively:
\[
  f_n (r) = 1 + O(e^{-m_\lam r} )
  {\hbox{ \quad and \quad }}
  a_n (r) = 1 + O(e^{-r}) \ .
\]
At the origin, $f_n(r)$ vanishes like $r^n$ and
$a_n(r)$ like $r^2$.
Hence $1-f_n(r)$ and $1-a_n(r)$ are well localized near the
origin.

The pair $(\psi^{(n)}, A^{(n)})$ is called the $n$-{\it vortex}
({\it magnetic} or {\it Abrikosov} (\cite{a,no}) in the case of
superconductors, and {\it Nielsen-Olesen} or {\it Nambu string} in
the particle physics case).  Note that $deg(\psi^{(n)}) = n$. No
other static solutions of the Ginzburg-Landau equations are
rigorously known, though there is a physical argument and
experimental evidence for the existence of vortex lattices -- the
Abrikosov lattices.

Observe that (in the present scaling)
the length scale for the magnetic field
(the {\em penetration depth}) is $1$, and the length
scale for the order parameter
(the {\em coherence length}) is $\frac{1}{m_{\lam}}$,
where $m_{\lam} = \min(\sqrt{2\lam},2)$.
More precisely, the following asymptotics for the
field components of the $n$-vortex
were established in~\cite{p} (see also~\cite{jt}):
as $r := |x| \ra \infty$,
\be
\label{eq:decay}
  \ba{c}
  j^{(n)}(x) = n \beta_{n} K_1(|x|)[1 + o(e^{-m_{\lam}r})] J\hat{x} \\
  B^{(n)}(r) = n \beta_{n} K_1(r)[1 - \frac{1}{2r} + O(1/r^2)] \\
  |1-f_n(r)| \leq c e^{-m_{\lam} r} \\
  |f_n'(r)| \leq c e^{-m_{\lam} r}.
  \ea
\end{equation}
Here $j^{(n)} := Im(\overline{\psi^{(n)}} \nabla_{A^{(n)}} \psi^{(n)})$
is the $n$-vortex supercurrent, and
$\beta_n > 0$ is a constant.
$K_1$ is the modified Bessel function of order $1$
of the second kind.  Since $K_1(r)$ behaves like
$c e^{-r} / \sqrt{r}$ for large $r$, we see
that the length scale for $j$ and $B$ is $1$.
Note that the two length length scales
$1/m_{\lam}$ and $1$ coincide at $\lam=1/2$.
Superconductors are referred to as {\it Type I}
if $\lam < 1/2$, and {\it Type II} if $\lam > 1/2$.

Consider test functions describing several vortices, with the
centers at points $z_1$, $z_2,\dots$ and with the degrees $n_1$,
$n_2$, $\ldots$ , glued together. An example of such a function
can be easily constructed as $v_{\zb,\chi} = (\psi_{\zb,\chi},
A_{\zb,\chi})$ with
\begin{equation}
  \psi_{\zb,\chi}(x) = e^{i\chi(x)} \prod_{j=1}^m
  \psi^{\nj}(x-z_j)
\end{equation}
and
\begin{equation}
  A_{\zb,\chi}(x) = \sum_{j=1}^m A^{\nj} (x-z_j) + \nabla \chi(x) \ ,
\end{equation}
where $\zb = (z_1,z_2,\dots)$ and $\chi$ is an arbitrary real
function yielding the gauge transformation (the integer degrees of
the vortices, $\underbar{n} = (n_1,\ldots,n_m)$, are suppressed in
the notation). Define the inter-vortex separation
\[
  R(\zb) := \min_{j \not= k} |a_j-a_k|.
\]
Since vortices are exponentially localized,
for large separation $R(\zb)$ (compared with
$[\min(m_\lam, 1)]^{-1}$) such test functions are
approximate -- but not exact --
solutions of the stationary Ginzburg-Landau equations.

When $\lam > 1/2$, we take $n_j = \pm 1$, since
vortices with $|n| \geq 2$ are known to be unstable
(\cite{gs}).

Now consider a time-dependent Ginzburg-Landau equation
with an initial condition $\vzo$ and ask
the following questions:  does the solution at time $t$
describe well-localized vortices at some locations
$\zb = \zb(t)$ (and with a gauge transformation $\chi = \chi(t)$)
and, if it does, what is the dynamic law of the vortex
centers $\zb(t)$ (and of $\chi (t)$)?

We describe here answers to these questions
for the superconductor model~(\ref{eq:gf}) and
Higgs model~(\ref{eq:mh}).  Precise statements
(Theorems~\ref{thm:gf} and~\ref{thm:mh})
are given in Section~\ref{sec:res}.

Consider the superconductor model~(\ref{eq:gf})
with initial data $(\psi_0,A_0)$
close to some $\vzo$ with
$e^{-R(\zb_0)}/\sqrt{R(\zb_0)} < \e$.
We show that the solution can be written as
\be
\la{eq:close}
  (\psi(t),A(t)) = \vzt + O(\e \log^{1/4}(1/\e))
\end{equation}
and that the vortex dynamics is governed by the system
\be
\la{eq:seff}
  \g_{n_j} \dot{z}_j = -\nabla_{z_j} W(\zb) + O(\e^2 \log^{3/4}(1/\e)).
\end{equation}
Here $\dot{z}_j$ denotes $dz_j/dt$,
$W(\zb) := \E_{GL}(\vz) - \sum_{j=1}^m E^{(n_j)}$,
where $E^{(n)} := \E_{GL}(\psi^{(n)},A^{(n)})$,
is the effective energy, and $\g_n$ are the numbers given by
\be
\la{eq:gamma'}
  \g_n := \frac{1}{2} \| \nabla_{A^{(n)}} \psi^{(n)} \|_2^2
         + \| curl A^{(n)} \|_2^2.
\end{equation}
In general, these statements hold only as long as the path
$\zb(t)$ does not violate a condition of large separation:
$R(\zb(t)) > \log(1/\e) + c$. In the {\it repulsive} case, when
$\lam > 1/2$ and $n_j = +1$ (or $n_j = -1$) for all $j$, the above
statements hold for all time $t$.  A precise statement is given in
Theorem~\ref{thm:gf}.

The leading-order term in the r.h.s of~(\ref{eq:seff}) is of order
$\e$ (see Lemma~\ref{prop:force} and Remark~\ref{rem:type1}). For
$\lam > 1/2$, the leading order of $W(\zb)$
for large $R(\zb)$ is:
\[
  W(\zb) \sim \sum_{j \not= k} (const) n_j n_k
  \frac{e^{-|z_j-z_k|}}{\sqrt{|z_j-z_k|}}
\]
(see Section~\ref{sec:reden}).

For the Higgs model equations
with initial data $(\psi_0,A_0)$
close to some $\vzo$
(and with appropriate initial momenta),
we show that
\be
\la{eq:close2}
  \|(\psi(t),A(t)) - \vzt\|_{H^1} +
  \|(\p_t \psi(t),\p_t A(t)) - \p_t \vzt\|_{L^2}
  = o(\sqrt{\e})
\end{equation}
with
\be
\la{eq:heff}
  \gamma_{n_j}\ddot{z}_j = -\nabla_{z_j} W(\zb(t)) + o(\e)
\end{equation}
for times up to (approximately) order
$\frac{1}{\sqrt{\e}}\log\left(\frac{1}{\e}\right)$.
Here $\ddot{z}_j(t)$ denotes $d^2 z_j(t)/dt^2$.
This result is stated precisely in Section~\ref{sec:res}
(see Theorem~\ref{thm:mh}).

The resulting dynamics of vortices induced by the field dynamics
of $(\psi, A)$ is called the {\it effective dynamics}.

We now outline some previous works on vortex dynamics, including
related works on the Gross-Pitaevski (or nonlinear Schr\"odinger)
equation \be \la{eq:gp}
  i \frac{\p \psi}{\p t}
  = -\Delta \psi + \frac{1}{\e^2}(|\psi|^2-1)\psi
\end{equation}
in a bounded domain,
used in the theory of superfluids
(see~\cite{tt}).  It is obtained from~(\ref{eq:schmax})
by setting $\g=i$ and $A=0$.
The landmark previous developments are summarized in the table
below

\begin{center}
\begin{tabular}{|l|l|l|l|} \hline
\hfill Type of&Superfluid&Superconductor&Higgs\\
\hfill Eqns&&&\\
\cline{1-1}
Type of\hfill &&& \\
Results\hfill &&&\\ \hline
Nonrigorous\qquad &Onsager `49 & Perez-Rubinstein &Manton `82\\
             &           & '83 $(\lambda\gg\frac12 )$
        & $(\lambda \approx \frac12 )$\\
             &           & E '84 $(\lambda\gg\frac12 )$&\\
\hline
Rigorous    &Colliander-&Demoulini-& Stuart `94\\
            &Jerrard '00&Stuart '97& $(\lambda\approx\frac12)$\\
            &F.-H. Lin-Xin `00 &$(\lambda\approx\frac12)$&\\ \hline
\end{tabular}
\end{center}

In more detail, non-rigorous results for the Ginzburg-Landau
equation (\ref{eq:gp}) without the magnetic component,
were obtained by
L. Onsager (\cite{on}),
A. Fetter (\cite{fe}),
R. Creswick and M. Morrison (\cite{cm}),
J. Neu (\cite{ne}),
L.M. Pismen and D. Rodriguez (\cite{pirod}),
D. Rodriguez, L.M. Pismen and L. Sirovich (\cite{rps}),
L.M. Pismen and J. Rubinstein (\cite{pirub}),
N. Ercolani and R. Montgomery (\cite{em}),
W. E (\cite{e}),
Yu. Ovchinnikov and I.M. Sigal (\cite{os}).

Rigorous results are contained in
J.E. Colliander and R.L. Jerrard (\cite{cj}),
F.-H. Lin and J. Xin (\cite{lx}),
based on Bethuel, Br\'ezis and H\'elein (\cite{bbh}).
Let $\psi^\epsilon$ be the solution of Eqn~(\ref{eq:gp})
with a ``low energy'' initial condition.
Then these papers show that as $\epsilon \to 0$,
the ``renormalized'' energy density
$$
  \frac{1}{|\log \e|} \left( \frac{1}{2} |\nabla \psi^{\e}|^2
  + \frac{1}{4\e^2}(|\psi^{\e}|^2-1)^2 \right)
$$
converges weakly to a sum of $\delta$-functions
located at points
$\zb(t) := \big( z_1 (t) , \dots , z_k (t)\big)$
which solve the Hamiltonian equation
$\dot{z} = J \nabla H(z)$
with appropriate initial conditions and Hamiltonian $H$.
Also \cite{cj} prove the
Bethuel-Br\'ezis-H\'elein type result
$\forall\rho>0$, as $\epsilon\to0 $
$$
  \min\limits_{\alpha\in [0,2\pi]}||\psi^\epsilon
  -e^{i\alpha} H_{\zb(t)}||_{H^1(T^2_\rho)}\to 0
$$
where $H_{\zb}$ is the Bethuel-Br\'ezis-H\'elein canonical
harmonic map with singularities at $z_1, \ldots, z_N$ and
$T^2_\rho=T^2 / \cup_i B_\rho(z_i)$,
and \cite{lx} show that the rescaled linear momentum
$Im(\bar{\psi}_\e \nabla \psi_\e)$
converges (on the time-scale $O(1)$) to a
solution of an incompressible Euler equation.
The results above describe the dynamics of the
vortex centers, but say nothing about the vortex
structure of the solutions.

In the magnetic case non-rigorous results were obtained in
N. Manton (\cite{m}) $(\lambda \approx \frac12)$,
M. Atiyah and N. Hitchin (\cite{ah}) $(\lambda\approx \frac12)$,
L. Perez and J. Rubinstein (\cite{pr}),
and W.E (\cite{e}).

Rigorous results were obtained in
D. Stuart (\cite{s}) $(\lambda\approx \frac12)$, and
S. Demoulini and D. Stuart (\cite{ds})
$(\lambda\approx \frac12)$.

Finally, we mention the recent results
\cite{ew,iww,abf,af,bf,cc,des,pego,rsk,
sw1,sw2,sw3,bp,bs,bj,fty,ty1,ty2,ty3,fgjs}
on interface, bubble, spike, and soliton dynamics.

The rest of the paper is organized as follows.
Ginzburg-Landau preliminaries are given
in Sections~\ref{sec:pre} and~\ref{sec:multi}.
The effective dynamics results described above
(Theorems~\ref{thm:gf} and~\ref{thm:mh})
are stated precisely in Section~\ref{sec:res}.
Theorem~\ref{thm:mh} is proved in Section~\ref{sec:mh},
and Theorem~\ref{thm:gf} in Section~\ref{sec:gf}.
The key technical estimates used in the proofs
are themselves proved in Section~\ref{sec:prop}.
Technical complications are relegated to
appendices (Sections~\ref{sec:UV}-~\ref{sec:app3}).

\bigskip

\noindent {\bf Notation.}
Here, and in what follows, $H^s$ denotes
the Sobolev space $H^s(\R^2 ; \C \times \R^2)$ (same for $L^2$,
etc.).
For $\mu = (\phi, \al), \nu=(\chi,\beta) \in L^2$,
$\lan \mu, \nu \ran$ denotes the real
$L^2$-inner product
\be
\la{eq:inner1}
  \lan \mu, \nu \ran := \int_{\R^2} \{
  Re (\overline{\mu} \nu) + \al \cdot \beta \}.
\end{equation}
Moreover, we will use the same symbol to denote the
real inner-product in $L^2 \times L^2$:
for $\xi = (\xi_1,\xi_2)$, $\eta = (\eta_1,\eta_2)$,
we write
\be
\la{eq:inner2}
  \lan \xi, \eta \ran :=
  \lan \xi_1, \eta_1 \ran + \lan \xi_2, \eta_2 \ran.
\end{equation}
$L^p$-norms are denoted with a subscript $p$:
$\| \cdot \|_p = \| \cdot \|_{L^p}$.
The letter $c$ will denote a generic
constant, independent of any small parameters present,
which may change from line to line.

%%%%%%%%%%%%%%%%%%%%%%%%%%%%%%%%%%%%%%%%%%%%%%%%%%
%%%%%%%%%%%%%%%    RESULTS    %%%%%%%%%%%%%%%%%%%%
%%%%%%%%%%%%%%%%%%%%%%%%%%%%%%%%%%%%%%%%%%%%%%%%%%

\section{Ginzburg-Landau preliminaries and results}

%----------------------------------------------------
\subsection{Ginzburg-Landau equations}
\la{sec:pre}

The Ginzburg-Landau energy functional $\E_{GL}$
(see~(\ref{eq:en})) is a smooth functional on the
following affine space of
configurations of degree $n$:
\[
  X^{(n)} := \{ (\psi,A) : \R^2 \ra \C \times \R^2 \;\; | \;\;
  (\psi,A)-(\psi^{(n)},A^{(n)}) \in H^1 \}
\]
where $(\psi^{(n)},A^{(n)})$
is the exact $n$-vortex solution of the
Ginzburg-Landau equations (see~(\ref{eq:en})).
The variational derivative $\E_{GL}'(\psi,A)$
is the (negative of the) right hand side
of the Ginzburg-Landau evolution
equations~(\ref{eq:gf}) (or~(\ref{eq:mh})).

With the notation $u = (\psi, A)$,
the superconductor model equations~(\ref{eq:gf})
can be written as
\[
  \p_t u(t) = -\E_{GL}'(u(t)).
\]
We consider solutions of~(\ref{eq:gf}) satisfying
$u=(\psi,A) \in C^1(\R^+ ; X^{(n)})$
(see~\cite{ds} for existence theory).

It is convenient to write
the Higgs model equations~(\ref{eq:mh}) as a
first-order Hamiltonian system.
Introduce the momenta
\[
  (\pi(t), E(t)) := (-\p_t \psi(t), -\p_t A(t))
\]
($E(t)$ is the electric field).  The Hamiltonian is \be
\la{eq:ham}
  \cH(\psi,A,\pi,E) := \E_{GL}(\psi,A)
  + \frac{1}{2}\int_{\R^2} \left\{ |\pi|^2 + |E|^2 \right\},
\end{equation}
a smooth functional on the space $X^{(n)} \times L^2$.
The space $X^{(n)} \times L^2$,
viewed as a real space, admits the
non-degenerate symplectic form
\be
\la{eq:sf}
  \omega(\xi, \eta) = \lan \xi, \J^{-1}\eta \ran
\end{equation}
where $\lan \cdot, \cdot \ran$ is the real inner product on the
tangent space to $X^{(n)} \times L^2$ defined
in~(\ref{eq:inner2}), and $\J$ is the symplectic operator
\[
  \J := \left( \begin{array}{cc}
  0 & {-\bfone} \\ {\bfone} & 0
  \end{array} \right)
\]
(in block notation).
Setting $w := (\psi,A,\pi,E)$,
the Higgs model~(\ref{eq:mh})
is equivalent to the equation
\be
\la{eq:hm}
  \p_t w(t) = \J \cH'(w(t)).
\end{equation}
We consider solutions in the space
$w \in C^1(\R ; X^{(n)} \times L^2)$
which conserve the Hamiltonian functional $\cH$
(see \cite{bm} for existence theory).

%----------------------------------------
\subsection{Multi-vortex configurations}
\label{sec:multi}

We begin by constructing a manifold of multi-vortex
configurations, made up of collections of widely-spaced vortices
``glued'' together. Such a collection is determined by $m \in
\Z^+$ vortex locations, $\zb = (z_1,\ldots,z_m) \in \R^{2m}$ and
$m$ vortex degrees, $\underbar{n} = (n_1,\ldots,n_m) \in \Z^m$,
associated with these locations (the latter will often be
suppressed in the notation), together with a gauge transformation.
So the manifold we construct may be parameterized by a subset of
$\R^{2m} \times$ $\{$gauge transformations$\}$.

Recall $(\psi^{(n)}, A^{(n)})$ denotes the equivariant, $n$-vortex
static solution of the Ginzburg-Landau equations
(see~(\ref{eq:vort})). To a triple $\zb \in \R^{2m}$,
$\underbar{n} \in \Z^m$, $\chi : \R^2 \ra \R$, we associate the
function \be \la{eq:vdef}
  \vz := (\psiz, \Az),
\end{equation}
where
\[
  \psiz(x) = e^{i\chi(x)} \prod_{j=1}^m \psi^{\nj}(x-z_j)
\]
and
\[
  \Az(x) = \sum_{j=1}^m A^{\nj}(x-z_j) +
  \nabla \chi(x).
\]
Here $(n_1, \ldots, n_m)$ are the fixed topological
degrees of the vortices, $n_j \in \Z \backslash \{0\}$.
For given $\zb \in \R^{2m}$, the gauge transformations
will be of the form
\[
  \chi(x) = \sum_{j=1}^m z_j \cdot A^{\nj}(x-z_j)
  + \chit(x)
\]
with $\chit \in H^2(\R^2;\R)$.
The gauge transformation is taken to be of this form
to ensure that $\vz$ lies in $X^{(n)}$.

Given a vortex configuration $\zb=(z_1,\ldots,z_m)$,
the inter-vortex distance is defined to be
\[
  R(\zb) := \min_{1 \leq j < k \leq m} |z_j-z_k|.
\]
To ensure that our multi-vortex configurations
are approximate solutions of the Ginzburg-Landau equations,
the inter-vortex separation will be taken large.

In the Higgs model case, momenta must be included.
To do this, we first introduce the
``almost zero-modes''.
Define the gauge ``almost zero-modes''
\be
\la{eq:modes1}
  \Gz_{\g} := \lan \g, \p_{\chi} \ran \vz
\end{equation}
for $\gamma : \R^2 \ra \R$,
and the gauge-invariant translational ``almost zero-modes''
\be
\la{eq:modes2}
  \Tz_{jk} :=
  (\partial_{z_{jk}} + \lan A^{\nj}_k(\cdot-z_j), \p_{\chi} \ran) \vz.
\end{equation}
From explicit expressions for $\Gz_\g$ and $\Tz_{jk}$
(see~(\ref{eq:texp}) and~(\ref{eq:gexp})), one can deduce that
$\Gz_\g, \Tz_{jk} \in H^s$,
provided $\g \in H^{s+1}$.
Then for momentum parameters
$\pb = (p_1,\ldots,p_m) \in \R^{2m}$ and
$\z \in H^1(\R^2;\R)$,
we define the $(\pi, E)$ (momentum)
component to be
\be
\la{eq:phidef}
  \phiz :=
  \sum_{j=1}^m p_j \cdot \Tz_j + \Gz_\z
  \in L^2.
\end{equation}
We will often denote the full set
of parameters by
$\s := (\zb,\chi,\pb,\z)$
and $\phiz$ by $\phi_\s$.

An important role will be played by the
interaction energy of a multi-vortex configuration
(see Section~\ref{sec:reden}):
\be
\la{eq:W}
  W(\zb) := \E_{GL}(\vz) - \sum_{j=1}^m E^{(n_j)}
\end{equation}
where, recall,
$E^{(n)} := \sum_{j=1}^m \E_{GL}(\psi^{(n)},A^{(n)})$.
Due to the gauge invariance of $\E_{GL}$,
this interaction energy is independent of the gauge
transformation $\chi$.

%-------------------------------------------------
\subsection{Main results}
\la{sec:res}

The main result in the superconductor model case
is as follows:
\begin{thm}
\label{thm:gf}
Suppose $\lam > 1/2$ and $n_j = +1$
(or $n_j=-1$) for $j=1,\ldots,m$.
There are $d_0,d,\e_0 > 0$ such that for
$0 < \e < \e_0$ the following holds:
let $(\psi(t),A(t))$ solve~(\ref{eq:gf}) with initial data
satisfying
\[
  \| (\psi_0,A_0) - \vzo \|_{H^1} < d_0 \e \log^{1/4}(1/\e)
\]
with $e^{-R(\zb_0)}/\sqrt{R(\zb_0)} < d_0 \e$.
Then for $t \geq 0$,
\[
  \| (\psi(t),A(t)) - \vzt \|_{H^1} < d \e \log^{1/4}(1/\e)
\]
for a path $\vzt \in M_{as}$ satisfying
\begin{equation}
\label{eq:gflaw}
  |\g_{n_j} \dot{z}_j(t) + \nabla_{z_j} W(\zb)| <
  d \e^2 \log^{3/4}(1/\e),
\end{equation}
\[
  \| \p_t \chi(t) - \sum_{j=1}^m \dot{z}_j(t)
  \cdot A^{\nj}(\cdot-z_j(t)) \|_{H^1} < d\e^2 \log^{3/4}(1/\e).
\]
\end{thm}
Here $\g_n$ is a positive constant,
given explicitly in~(\ref{eq:gamma'}).

For the Higgs model equations,
we have the following result:
\begin{thm}
\la{thm:mh}
Suppose $\lam > 1/2$ and $n_j = +1$
(or $n_j = -1$) for $j=1,\ldots,m$.
Let $\al(\e)$ be a function satisfying
$\sqrt{\e} < \al(\e) << \log^{-1/4}(1/\e)$.
There are $d_0,d,\tau, \e_0 > 0$ such that
for $0 < \e < \e_0$, the following
holds: let $w(t) = (\psi(t),A(t),\pi(t),E(t))$
solve~(\ref{eq:hm}), with initial data satisfying
\[
  \| (\psi_0,A_0) - \vzo \|_{H^1} +
  \| (\pi_0,E_0) - \phizo \|_2 < d_0 \e \log^{1/2}(1/\e)
\]
with
$e^{-R(\zb_0)}/\sqrt{R(\zb_0)} + |\pb_0|^2 + \|\zeta_0\|_2^2
< d_0 \e$.
Then for
$0 \leq t \leq \frac{\tau}{\sqrt{\e}} \log \left(
\frac{\al(\e)}{\sqrt{\e}} \right)$,
\begin{equation}
\label{eq:mhbound}
  \| (\psi(t),A(t)) - \vzt \|_{H^1} +
  \| (\pi(t),E(t)) - \phi_{\s(t)} \|_2 <
  d \al(\e) \sqrt{\e} \log^{1/2}(1/\e)
\end{equation}
for a path
$\s(t) = (\zb(t),\chi(t),\pb(t),\z(t))$
satisfying, for all $j$,
\be
\la{eq:mhlaw}
\begin{split}
  &|\dot{z}_j - p_j| +
  |\g_{n_j} \dot{p}_j +  \nabla_{z_j} W(\zb)| <
  \e \al(\e) \log^{1/2}(1/\e) = o(\e) \\
  &\| \p_t \chi - \sum_{j=1}^m \dot{z}_j \cdot A^{\nj}(x-z_j) -\z \|_{H^1}
  + \| \p_t \z \|_{H^{1-s}} < \e \al(\e) \log^{1/2}(1/\e) = o(\e)
\end{split}
\end{equation}
for any $s>0$.
\end{thm}
\begin{rem}
The inter-vortex force is of size $\e$:
$\nabla W(\zb) = O(\e)$
(see Lemma~\ref{prop:force} and Remark~\ref{rem:type1}).
\end{rem}
\begin{rem}
The condition $\lam > 1/2$ and $n_j=+1$ in Theorems~\ref{thm:gf}
and~\ref{thm:mh} ensures that the inter-vortex interaction
is repulsive, and therefore that the
inter-vortex separation does not become too small in the given time interval.
In fact the theorems apply, without these restrictions,
for any initial vortex configuration
whose evolution (namely~(\ref{eq:mhlaw}) or~(\ref{eq:gflaw}))
preserves an appropriate large-separation condition.
In the Type II case ($\lambda > 1/2$), this condition is
$e^{-R(\zb(t))}/\sqrt{R(\zb(t)} < \e$
(and $|\pb(t)| + \|\z(t)\|_{L^2} < \e$ in the
Higgs model case).
In the Type-I case ($\lam < 1/2$), this condition
must be appropriately modified
(see Remark~\ref{rem:type1} of Section~\ref{sec:reden}).
\end{rem}
%\begin{rem}
%One can obtain a longer time interval
%$(O(1/\e))$ and a smaller bound on
%$(\psi,A)-\vz$ (in fact $O(\e)$)
%for a restricted class of initial
%vortex configurations.
%This class is characterized as follows:
%under the leading order evolution
%of~(\ref{eq:mhlaw}), vortex velocities
%remain $O(\e)$.
%\end{rem}
%\begin{rem}
%The size $d_0 \e$ of the initial perturbation,
%and the small parameter
%$R(\zb_0) e^{-R(\zb_0)}$
%do not necessarily have to be the same.
%\end{rem}
\begin{rem}
In Theorem~\ref{thm:mh},
since $|\dot{z}_j| \leq c \sqrt{\e}$ over the time interval
$0 \leq t \leq T = \frac{\tau}{\sqrt{\e}} \log
\left( \frac{\al}{\sqrt{\e}} \right)$,
vortices can move a distance
\[
  \sqrt{\e} T = \ta
  \log{\frac{\al(\e)}{\sqrt{\e}}} \sim R(\zb(0)) >> 1.
\]
\end{rem}

%%%%%%%%%%%%%%%%%%%%%%%%%%%%%%%%%%%%%%%%%%%%%%%%%%%%%%%
%%%%%%%%    Proofs of vortex motion laws      %%%%%%%%%
%%%%%%%%%%%%%%%%%%%%%%%%%%%%%%%%%%%%%%%%%%%%%%%%%%%%%%%

\section{Proofs}
\la{sec:proofs}

We start by proving Theorem~\ref{thm:mh} for the
Higgs model in Section~\ref{sec:mh}.  The proof
is considerably more involved than that of Theorem~\ref{thm:gf}
for the superconductor model.  The latter proof is
sketched in Section~\ref{sec:gf}.

The proof of Theorem~\ref{thm:mh} given in the following section
is based on a series of propositions and lemmas.
Propositions~\ref{prop:proj1}-\ref{prop:tub} summarize our
geometric construction, and Lemmas~\ref{lem:nl1}-\ref{lem:repulse},
whose proofs are left to Section~\ref{sec:prop}, provide the
(elementary) analytic building blocks. Several technical lemmas are
relegated to appendices.

%-------------------------------------------------------
\subsection{Effective dynamics of vortices: Higgs model}
\la{sec:mh}

In this section we prove Theorem~\ref{thm:mh}.
Let $w(t)$ solve~(\ref{eq:hm}) with
$w \in C^1(\R;X^{(n)} \times L^2)$.
In what follows, we denote
\[
  X := H^1(\R^2;\C \times \R^2) \times L^2(\R^2;\C \times \R^2).
\]

%-----------------------------------------------------
\subsubsection{Manifold of multi-vortex configurations}

We begin by defining the manifold of multi-vortex configurations. Let
\[
  \Si := \{ (\zb,\chi,\pb,\z)
  \; | \; \zb \in \R^{2m}, \chi - \zb \cdot A_{\zb}
  \in H^2(\R^2;\R), \pb \in \R^{2m},\z \in H^1(\R^2;\R)\},
\]
where $\zb \cdot A_{\zb} := \sum_{j=1}^m z_j \cdot A_j$ with
$A_j(x) := A^{\nj}(x-z_j)$. This set is a manifold
under the explicit parametrization map $\del: Y_{2,1} \ra \Si$
defined by
\be
\la{eq:deldef}
  \del : (\zb,\chit,\pb,\z) \mapsto
  (\zb, \chit + \zb \cdot A_{\zb}, \pb, \z).
\end{equation}
Here
\[
  Y_{r,s} := \R^{2m}
  \times H^{r}(\R^2;\R)\times \R^{2m}
  \times H^{s}(\R^2;\R).
\]
We define an open domain in $\Si$ by
\[
 \Si_\e := \{ (\zb,\chi,\pb,\z) \in \Si
  \; | \; e^{-R(\zb)}/\sqrt{R(\zb)} < \e, \;
  |\pb| + \|\z\|_{H^1} < \sqrt{\e} \}.
\]

For each $\s := (\zb,\chi,\pb,\z) \in \Si$,
introduce the multi-vortex configuration
\be
\la{eq:wdef}
  w_\s := ( \vz, \phi_\s) \in
  X^{(n)} \times L^2
\end{equation}
(recall $\vz$ and $\phi_\s$ are defined in~(\ref{eq:vdef})
and ~(\ref{eq:phidef})).  Finally, we define the space
\[
  M_{mv} := \{ w_\s \; | \; \s \in \Si_\e \}
  \subset X^{(n)} \times L^2.
\]
The map $\g : \Si_\e  \ra X^{(n)} \times L^2$
given by $\g : \s  \to w_\s$,
parameterizes $M_{mv}$, so that $M_{mv} = \g(\Si_\e)$.
It is easy to check that $\g$ is $C^1$.
It is shown in Section~\ref{sec:UV} that for all
$\s \in \Si_\e$, its Fr\'echet derivative
$D \g(\s) : T_\s \Si_\e \to X$ is one-to-one.
Hence $M_{mv}$ is a manifold.
% ?? based on $Y_{2,1}$ ??

For each $\s \in \Si_\e$, the tangent space to
$M_{mv}$ at $w_\s$ will be denoted by $T_{w_\s} M_{mv}$.
It can be identified with a subspace of $X$; specifically,
$T_{w_\s} M_{mv} = D \g(\s) (T_\s \Si_\e)$.

For use in computations and estimates below,
we introduce convenient bases in
$T_\s \Si_\e$ and $T_{w_\s}M_{mv}$.
In terms of the coordinates in~(\ref{eq:deldef}),
the basis in $T_\s \Si_\e$ is given by
\be
\la{eq:Ybasis}
  \{\p_{z_{ij}} + \lan z_i \cdot \p_{x_j} A_i, \p_{\chit} \ran,
  \p_{\chit(x)}, \p_{p_{ij}}, \p_{\z(x)} \}.
\end{equation}

We denote the coordinates of
$\s' \in T_\s \Si_\e$ in this basis
by $\s'_{coord} \in Y_{2,1}$.
Define the map $\G_\s : Y_{2,1} \to X$ by
\be
\la{eq:Gamma1}
  \G_\s \s'_{coord} := D \g(\s) \s'.
\end{equation}
For $\s(t)$ a path in $\Si_\e$,
this definition implies
\[
  \p_t w_{\s(t)} = \G_{\s(t)} \dot{\s}(t),
\]
where $\dot{\s}(t)$ is the coordinate
representation of the vector
$\p_t \s(t) \in T_{\s(t)} \Si_\e$:
\[
  \dot\s(t) := ({\dot{\zb}(t)},{\p_t^{\zb(t)}
  \chi(t)},{\dot{\pb}(t)},{\p_t \z(t)}),
\]
with $\dot{\zb}(t) := d \zb(t)/dt$,
$\dot{\pb}(t) := d \pb(t)/dt$, and
\[
  \p^{\zb(t)}_t \chi(x,t) :=
  \p_t \chi(x,t) - \sum_{j=1}^m
  \dot{z}_j(t) \cdot A^{\nj}(x-z_j(t)).
\]

Let $\p_{z_{ij}}^A := \p_{z_{ij}} + \lan A_{ij}, \p_{\chi} \ran.$
The basis for $T_{w_\s} M_{mv}$ (which is the image of
the basis~(\ref{eq:Ybasis}) under $D\g(\s)$) is given by:
\be
\la{eq:tan}
  \ta^z_{ij} := \p_{z_{ij}}^A w_\s, \quad
  \ta^p_{ij} := \p_{p_{ij}} w_\s, \quad
  \ta^\chi_x := \p_{\chi(x)} w_\s, \quad
  \ta^\z_x := \p_{\z(x)} w_\s.
\end{equation}
Note that the tangent vector $\ta^z_{ij}$ is defined by
differentiating $w_\s$ {\it covariantly}.
The point here is that $(\p_z^A)^m w_\s$
lies in $H^1 \times L^2$ for any $m$,
while $\p_z w_\s$ does not.
Explicit expressions for these tangent vectors
are given in~(\ref{eq:tan1'})-~(\ref{eq:tan4'}).

For a vector $\al \in \R^{2m}$, we will set
$\al \cdot \ta^\# := \sum_{ij} \al_{ij} \ta^\#_{ij}$ for
$\ta^\#_{ij} = \ta^z_{ij}, \ta^p_{ij}$, and for a function $\g$,
set $\lan \g, \ta^\# \ran := \int \g(x) \ta^\#_x dx$ for
$\ta^\#_x = \ta^\chi_x, \ta^\z_x$.
As a result of these definitions, and the relation
\be
\la{eq:dt}
\begin{split}
   \p_t &= \dot{\zb} \cdot \p_{\zb}
   + \lan \p_t\chi, \p_\chi \ran
   + \dot{\pb} \cdot \p_{\pb}
   + \lan \p_t\z, \p_\z \ran \\
   &= \dot{\zb} \cdot \p_{\zb}^A
   + <\p^{\zb(t)}_t \chi, \p_\chi>
   + \dot{\pb} \cdot \p_{\pb}
   + \lan \p_t\z, \p_\z \ran,
\end{split}
\end{equation}
we have
\be
\la{eq:Gamma2}
  \G_\s \s'_{coord} = \zb' \cdot \ta^z + \lan  \chi', \ta^\chi \ran
  + \pb' \cdot \ta^p + \lan \z', \ta^\z \ran,
\end{equation}
where $\s'_{coord} = (\zb', \chi',\pb',\z') \in Y_{2,1}$.

In what follows, all of our computations are done
in these bases, and we omit the subscript
``coord'' from the coordinate representation
$\s'_{coord}$ of a vector $\s' \in T_\s \Si_\e$.

%--------------------------------------------------
\subsubsection{Reduced (vortex) Hamiltonian system}

As was discussed above, the Maxwell-Higgs equations constitute a
Hamiltonian system on the phase-space $X^{(n)} \times L^2$ with
Hamiltonian ~(\ref{eq:ham}). Our goal below is to project this
Hamiltonian system onto the manifold $M_{mv}$ (more precisely,
onto $TM_{mv}$) with the smallest error possible. Below we
describe an equivalent Hamiltonian structure on the parameter
space $Y_{2,1}$ which is used in our analysis.
We begin by setting
\[
  X_{r,s} := H^r(\R^2; \C \times \R^2) \times
             H^s(\R^2; \C \times \R^2)
\]
(note that $X = X_{1,0}$).
The operator $\G_\s$ has adjoint $\Lam_\s$ (with respect
to the $\R^{2m} \times L^2 \times \R^{2m} \times L^2$
inner-product on $Y_{r,s}$, and the real $L^2 \times L^2$
inner-product on $X_{r,s}$) given by
\be
\la{eq:lamdef}
  \Lam_\s : \xi \mapsto
  \lan D_\s w_\s, \xi \ran
\end{equation}
or, in our coordinates in $T_\s\Si_{\e}$,
\be
\la{eq:ladef'}
  \Lam_\s : \xi \mapsto
  (\lan \ta^z_{ij}, \xi \ran, \lan \ta^\chi_x, \xi \ran, \lan
  \ta^p_{ij}, \xi \ran, \lan \ta^\z_x, \xi \ran).
\end{equation}
It is shown in Section~\ref{sec:UV} that $\G_\s$ and $\Lam_\s$
are bounded uniformly in $\s \in \Si_\e$ between the following spaces:
\be
\la{eq:gamap}
  \G_\s : Y_{r,s} \ra X_{r-1,s-1},
\end{equation}
for any $r$ and $s$ satisfying $min (r,1) > s-1$, and
\be
\la{eq:lamap}
  \Lam_\s : X_{r,s} \ra Y_{r-1,s-1},
\end{equation}
for any $r$ and $s$ satisfying $min(s,1) > r-1$. (In the
Physics literature the operators $\G_\s$ and $\Lam_\s$ are called
bra and ket vectors, with notation $\G_\s = |D w_\s \ran$ and
$\Lam_\s = \lan D w_\s|$.)

Define the operators
\be
\la{eq:UVdef}
   V_\s := \Lam_\s \J^{-1} \G_\s : Y_{r,s} \ra Y_{s-2,r-2},
\end{equation}
where $s-1 < min(r,2)$.

Relation~(\ref{eq:UVdef}) shows that $V_\s^* = -V_\s$ in the sense
of the $L^2$ inner product. The operators $V_\s$ define a
symplectic form on $Y_{2,1}$ by
\[
  \omega_{red}(\s',\s'')(\s) :=
  \lan \s', V_\s \s'' \ran.
\]
The non-degeneracy of this symplectic form follows from:
\begin{prop}[non-degeneracy of reduction]
\label{prop:proj1}
For $\e$ sufficiently small, and $\s \in
\Si_\e$, the operator $V_\s$ is invertible.
\end{prop}
{\bf Proof}:
The invertibility of the operator $V_\s$ for sufficiently small
$\e$  follows from the following expression,
shown in Section~\ref{sec:UV}:
\be
\la{eq:Vs}
  V_\s =
  \left( \begin{array}{cc}
  R_1 & -B \\ B^* & R_2
  \end{array} \right)
\end{equation}
where
\[
  B = \left( \begin{array}{cc} D & O(\e\log^{1/2}(1/\e)) \\
  O(\e\log^{1/2}(1/\e)) & K \end{array} \right), \quad
  R_1 = \left( \begin{array}{cc} 0 & O(\sqrt{\e}) \\
  -O(\sqrt{\e})^* & 0 \end{array} \right),
\]
\[
  R_2 = \left( \begin{array}{cc} 0 & O(\e\log^{1/2}(1/\e)) \\
  -O(\e\log^{1/2}(1/\e))^* & 0 \end{array} \right),
\]
$O(\sqrt{\e})$ stands for an operator whose norm is bounded
by $c \sqrt{\e}$, $D$ is a matrix of the form
$D_{jk,lm} = \g_{n_j} \del_{jl} \del_{km} + O(\e\log^{1/2}(1/\e))$,
and $K$ is the operator $K:= -\Delta + |\psiz|^2$.
%Here $\del = \del(\s) := |\pb| + \|\z\|_{H^1}$.
%Recall that $\del < \sqrt{\e}$ if $\s \in \Si_\e$.
Since $D$ and $K$ are invertible, the operators
$U_\s$ and $V_\s$ are obviously
invertible if $\e$ is sufficiently small.
$\Box$

The symplectic form $\omega_{red}(\s',\s'')$ and the reduced (vortex)
Hamiltonian $h(\s) := \cH(w_\s)$ give a
reduced Hamiltonian system on $Y_{2,1}$.
The corresponding Hamiltonian equation is
\[
  \dot \s = V_\s^{-1} Dh(\s).
\]
This equation will turn out to be the leading-order equation for the
dynamics of the parameters $\s$. The next proposition computes the
Hamiltonian $h(\s)$ explicitly.
\begin{prop}
\la{prop:redham}
If $e^{-R(\zb)}/\sqrt{R(\zb)} < \e$, then
\be
\la{eq:expansion}
\begin{split}
  h(\s) := \cH(w_\s) &= \sum_{j=1}^m E^{(n_j)}
  + W(\zb) + \frac{1}{2} \sum_{j=1}^m \g_{n_j} |p_j|^2
  + \frac{1}{2} \lan \zeta, K \zeta \ran \\
  &\quad + O(\e\log^{1/2}(1/\e)(|p|^2 + \|\z\|_2^2))
\end{split}
\end{equation}
where, recall, $E^{(n)} := \E_{GL}(\psi^{(n)},A^{(n)})$,
$\g_n = \frac{1}{2} \| \nabla_{A^{(n)}} \psi^{(n)} \|_2^2
+ \| curl A^{(n)} \|_2^2$,
and $W(\zb)$ is defined in~(\ref{eq:W})
(it is computed to leading order in Section~\ref{sec:reden}).
Further, we have
\be
\la{eq:dexpansion}
  Dh(\s) := D_\s\cH(w_\s) =
  (\nabla W(\zb), 0, \underline{\g} \cdot \pb, K \zeta)
  + O(\e\log^{1/2}(1/\e)(|p| + \|\z\|_2)),
\end{equation}
where $\underline{\g} \cdot \pb$ denotes
$(\g_1 p_1,\ldots,\g_m p_m)$.
\end{prop}
{\it Proof}. We begin with auxiliary computations establishing the
approximate orthogonality of the tangent vectors introduced above.
To this end, we record explicit expressions for $\Tz_{jk}$ and
$\Gz_\g$, which follow readily from definitions~(\ref{eq:modes1})
and~(\ref{eq:modes2}):
\be
\la{eq:texp}
  \Tz_{jk} = -(e^{i\chi}
  [\prod_{l \not= j} \psi^{\nl}(x-z_l)]
  (\nabla_{A_k} \psi)^{\nj}(x-z_j), \;
  B^{\nj}(x-z_j) e_k^{\perp})
\end{equation}
and
\be
\la{eq:gexp}
  \Gz_\g = (i \g \psiz, \; \nabla \g).
\end{equation}
Here $B^{(n)} = \nabla \times A^{(n)}$ is the $n$-vortex magnetic
field, $e_1^{\perp}:=(0,1)$ and $e_2^{\perp}:=(-1,0)$.

By the above explicit expressions, the exponential decay
estimates~(\ref{eq:decay}), and Lemma~\ref{lem:*}, we see
\[
  |\lan \Tz_{jr}, \Tz_{ks} \ran| \leq c\e \log^{1/2}(1/\e)
\]
when $j \not= k$.  When $j=k$, we compute
\[
\begin{split}
  \lan\Tz_{jr}, \Tz_{js}\ran &=
  \lan({\dA}_r\psi)_j, ({\dA}_s\psi)_j\ran  \\
  &\quad + \lan B_j J\hat{e}_r, B_j J\hat{e}_s \ran + O(\e),
\end{split}
\]
and the leading term is easily computed to be $\g_{n_j} \del_{rs}$
where $\g_{n_j}$ is given in ~(\ref{eq:gamma'}). Thus we have the
approximate orthogonality relation
\be
\la{eq:Dexp}
  D_{jk,lm} := \lan \Tz_{jk}, \Tz_{lm} \ran
  = \g_{n_j} \del_{jl} \del_{km} + O(\e\log^{1/2}(1/\e)).
\end{equation}
A similar computation yields
\be
\la{eq:mixorth}
  | \lan\Tz_{jk}, \Gz_\g \ran | \leq c \e \log^{1/2}(1/\e)
  \| \g \|_2.
\end{equation}
Finally, the corresponding relation for the approximate gauge
modes (see~(\ref{eq:modes1})) is \be \la{eq:gorth}
  \lan \Gz_\g, \Gz_\z \ran = \lan \g, (-\Delta + |\psiz|^2) \z \ran,
\end{equation}
a straightforward calculation.

Now using $w_\s = (\vz, \phi_\s)$ (with $\phi_\s$ defined
in~(\ref{eq:phidef})), and
\[
  \| \phi_\s \|_2^2 = p_{jk} p_{rs} \lan \Tz_{jk}, \Tz_{rs} \ran
  + \lan \Gz_\z, \Gz_\z \ran
  - 2 p_{jk} \lan \Tz_{jk}, \Gz_\z \ran
\]
together with~(\ref{eq:Dexp})-(\ref{eq:gorth}), we
obtain~(\ref{eq:expansion}) and~(\ref{eq:dexpansion}).
% ?? do we get dexpansion ??
$\Box$

%---------------------------------------------------------
\subsubsection{Projections $Q_\s$}

Here we construct operators $Q_\s$ used to engineer a convenient
splitting of~(\ref{eq:hm}). We define the operator $Q_\s : X \ra
T_{w_\s} M_{mv}$ as
\be \la{eq:exprQ}
  Q_\s := \G_\s V_\s^{-1} \Lam_\s \J^{-1}.
\end{equation}
Due to the expression for $V_\s$ in~(\ref{eq:UVdef}),
we see that $Q_\s$ is a projection,
$Q_\s^2 = Q_\s$, and it satisfies
\be
\la{eq:qc}
  Ker Q_\s = ( \J T_{w_\s} M_{mv})^{\perp}
\end{equation}
and
\be
\la{eq:Q*}
  Q_\s^* = -\J Q_\s \J.
\end{equation}

Finally, we list two estimates which follow
readily from the definitions above:
\be
\la{eq:Qsob}
  \| Q_\s \|_{X_{r,s} \ra X_{r,s}}
  \leq c
\end{equation}
for any $r$ and $s$ satisfying $s < min(r+1,1)$,
% ?? further conditions on exponents ??
and for $\s = \s(t)$ a path in $\Si_\s$,
\be
\la{eq:Q'}
  \| [Q_\s, \p_t] \|_{X \to X}
   \leq c \| \dot \s \|_{Y_{1,0}}
\end{equation}
where, recall,
$\dot \s = (\dot{\zb}, \p_t^{\zb} \chi, \dot{\pb}, \p_t \z)$.
% ?? check it ??
To obtain this estimate one uses relation~(\ref{eq:dt}).

%------------------------------------------------------
\subsubsection{Splitting}

The next proposition establishes a coordinate system
(adapted to the projection $Q_\s$)
on a tubular neighbourhood of $M_{mv}$.
Let, for $0 < d_0 < 1$,
\[
  \Si^0_\e := \{ (\zb,\chi,\pb,\z) \in \Si
  \; | \; e^{-R(\zb)}/\sqrt{R(\zb)} < d_0\e,
  |\pb| + \|\z\|_{H^1} < d_0 \sqrt{\e} \}
\]
(which parameterizes a manifold somewhat
smaller than $M_{mv}$).  Set
\[
  U_\del := \{ w \in X^{(n)} \times L^2 \; | \;
  \|w - w_\s\|_X < \del,
  \mbox{ for some } \s \in \Si^0_\e \}.
\]
\begin{prop}[coordinates]
\la{prop:tub}
For $\e$ sufficiently small,
there is $\del >> \e$, and a $C^1$ map
\[
  S : U_\del \ra \Si_\e
\]
satisfying $Q_{S(w)}(w - w_{S(w)}) = 0$ for $w \in U_\del$.
Moreover, $DS(w)$ is bounded uniformly in $w \in U_\del$.
\end{prop}

\noindent
{\bf Proof}:
The proof is an application of the implicit
function theorem. Define
\[
  g : U_\del \times \Si_\e^0 \ra Y_{-1,0}
\]
by
\[
  g(w,\s) := \Lam_\s \J^{-1}(w-w_\s).
\]
One can check that this is a $C^1$ map. Obviously, $g(w_\s,\s) =
0$. Note that, due to~(\ref{eq:Gamma1}), $D_{\s}g(w_\s,\s) :
Y_{2,1} \ra Y_{-1,0}$ is given by
\[
  D_{\s} g(w_\s,\s) = - \Lam_\s \J^{-1}\G_\s = - V_\s.
\]
which is invertible for $\s \in \Si_\e$ with $\e$ sufficiently
small. So the implicit function theorem applies to provide a $C^1$
map $w \mapsto S(w)$ from an $H^1$-ball of size $\del$ of a given
$w_\s \in M_{mv}$ into $\Si$, satisfying $g(w,S(w)) = 0$.
Allowing $\s$ to vary in $\Si^0_\e$, we can construct such a ball
about any such $w_\s$.

Using the definitions~(\ref{eq:Gamma1}) and~(\ref{eq:lamdef}) of
the operators $\G_\s$ and $\Lam_\s$, the explicit
expressions~(\ref{eq:tan1'})-~(\ref{eq:Fdef}) for the basis
$\{\ta^z,\ta^\chi,\ta^p,\ta^\z\}$, and expression~(\ref{eq:Vs}),
one can check the following: there is $\del_0$ independent of
$\s \in \Si_\e$ such that for all $w \in B_X(w_\s;\del_0)$, the norms
\[
\begin{split}
  &\| V_\s^{-1} \|_{Y_{-1,0} \to Y_{2,1}}, \quad
  \| \G_\s \|_{Y_{2,1} \to X_{1,0}}, \quad
  \| \Lam_\s \|_{X_{0,1} \to Y_{-1,0}}, \\
  &\| D_\s \Lam_\s \|_{X_{0,1} \times Y_{2,1} \to Y_{-1,0}},
  \quad
  \| D_\s \G_\s \|_{Y_{2,1} \times Y_{2,1} \to X_{1,0}},
  \quad
  \| D_\s^2 \Lam_\s \|_{Y_{2,1} \times Y_{2,1}
    \times X_{0,1} \to Y_{-1,0}}
\end{split}
\]
are bounded uniformly in $\s$. This fact implies that the balls
on which the maps $S$ are defined can be taken to be of uniform size
$\delta << \sqrt{\e}$, which implies
$S(w) \in \Si_\e$. Thus we obtain a well-defined $C^1$ map $S: w
\mapsto S(w)$ from the tubular neighborhood $U_{\del}$ into
$\Sigma_\e$, with $g(w,w_{S(w)}) = 0$. This map obviously
satisfies also $Q_{S(w)}(w - w_{S(w)}) = 0$. The uniform
boundedness of $DS(w)$ follows readily from the formula $DS(w) =
-[D_\s g(w,S(w))]^{-1} D_w g(w,S(w))$ and the uniform estimates
mentioned above. $\Box$

Now suppose $w(t)$ solves the Higgs model equations~(\ref{eq:hm})
with initial data $w(0)=w_0$ as specified in Theorem~\ref{thm:mh}.
In particular, we have $w(0) \in U_\del$.
Let $0 < T_1 \leq \infty$ be the time
of first exit of $w(t)$ from $U_\del$.
For $0 \leq t < T_1$ we may write
\be
\la{eq:splitter}
  w(t) = w_{\s(t)} + \xi(t)
\end{equation}
with $w_{\s(t)} \in M_{mv}$, and
$Q_{\s(t)} \xi(t) \equiv 0$
(by choosing $\s(t) = S(w(t))$).
By our choice of initial data,
\be
\la{eq:xi(0)a}
  \| \xi(0) \|_X < c \| \xi_0 \|_X < c d_0 \e \log^{1/4}(\e),
\end{equation}
where $\xi_0 := w(0) - w_{\s_0}$.
Indeed, using~(\ref{eq:splitter}) and the
equation $w(0) = w_{\s_0} + \xi_0$, we find
\be
\la{eq:xi(0)}
  \xi(0) = w_{\s_0} - w_{\s(0)} + \xi_0.
\end{equation}
Next, since $\s(0) = S(w(0))$ and $\s_0 = S(w_{\s_0})$
(see Proposition~\ref{prop:tub}), and since
$w(0)-w_{\s_0} = \xi_0$,
Proposition~\ref{prop:tub} gives
\[
  \| \s_0 - \s(0) \|_{Y_{2,1}} \leq c \| \xi_0 \|_X.
\]
The last estimate, together with the estimate
$\|D_\s w_\s \|_{Y_{2,1} \to X} \leq c$ implies that
\[
  \|w_{\s_0} - w_{\s(0)}\|_X \leq c \| \xi_0 \|_X,
\]
which, together with~(\ref{eq:xi(0)}),
yields~(\ref{eq:xi(0)a}).

%-------------------------------------------
\subsubsection{Effective dynamics}

Insert the decomposition~(\ref{eq:splitter})
into the equations~(\ref{eq:hm}) and expand in a
Taylor series to obtain
\be
\la{eq:exp}
  \p_t w_\s + \p_t \xi =
  \J[ \cH'(w_\s) + L_\s \xi + N_\s(\xi) ]
\end{equation}
where $L_\s := \cH''(w_\s)$ is the Hessian of $\cH$
at $w_\s$, and
\[
  N_\s(\xi) := \cH'(w_\s+\xi) - \cH'(w_\s) - L_\s \xi
\]
consists of the terms nonlinear in $\xi$.
Apply the projection $Q_\s$ to~(\ref{eq:exp})
and use $Q_\s \p_t w_\s = \p_t w_\s$
(since $\p_t w_\s \in T_{w_\s} M_{mv}$) to obtain
\be
\la{eq:ml}
  \p_t w_\s - Q_\s \J \cH'(w_\s)
  = Q_\s[\J L_\s \xi - \p_t \xi + \J N_\s(\xi)].
\end{equation}
This equation governs the effective dynamics of the
parameters $\s(t)$.  The terms of leading order
are on the left hand side.
We now show, starting with the nonlinear
term, that the right hand side is of lower order.
\begin{lem}[nonlinear estimate 1]
\la{lem:nl1}
For $\s \in \Sigma_\e$,
and $\xi := (\xi_1,\xi_2) \in H^1 \times L^2$,
\[
  N_\s(\xi) = \left( \begin{array}{c}
  (N_\s)_1(\xi_1) \\ 0 \end{array} \right)
\]
with
\[
  \| (N_\s)_1(\xi_1) \|_{H^{-s}}
  \leq c_s(\|\xi_1\|_{H^1}^2 + \|\xi_1\|_{H^1}^3)
\]
for any $s>0$.
\end{lem}
This lemma is proved in Section~\ref{sec:tay}.
From now on, we fix $s > 0$ ($s=1/2$, say).
Thus by~(\ref{eq:Qsob}), we have
\[
  \| Q_\s \J N_\s(\xi) \|_{L^2 \times H^{-s}} \leq
  c(\|\xi\|_X^2 + \|\xi\|_X^3).
\]

To minimize writing in the rest of this section,
we make the additional assumption
\be
\la{eq:xismall}
  \| \xi \|_X < 1,
\end{equation}
which we shall justify later.
Then the above estimate becomes
\be
\la{eq:inq1}
 \| Q_\s \J N_\s(\xi) \|_{L^2 \times H^{-s}} \leq
  c\|\xi\|_X^2.
\end{equation}

Using the fact that $Q_\s \xi \equiv 0$ and the
bound~(\ref{eq:Q'}), we have
\be
\la{eq:inq2}
  \|Q_\s \p_t \xi\|_X
  \leq c \| \dot \s \|_{Y_{1,0}}
  \|\xi\|_X.
\end{equation}

To bound the remaining term on the right hand side
of~(\ref{eq:ml}), we need the following lemma
whose proof is given in Section~\ref{sec:zero}.
\begin{lem}[approximate zero-modes]
\la{lem:zero}
For $\s \in \Si_\e$ and any
$\beta \in L^2 \times L^2$, we have
\be
\la{eq:approx1}
  \| L_\s Q_\s \beta \|_{L^2 \times L^2}
  \leq c \sqrt{\e} \| \beta \|_{L^2 \times L^2}.
\end{equation}
\end{lem}

Fix $\eta \in L^2 \times L^2$.
Using the symmetry of $L_\s$,
and~(\ref{eq:approx1}), we have
\[
\begin{split}
  | \lan \eta, Q_\s \J L_\s \xi \ran | &=
  | \lan L_\s Q_\s \J \eta, \xi \ran |
  \leq \|\xi\|_{H^1 \times L^2}
  \| L_\s Q_\s \J \eta \|_{H^{-1} \times L^2} \\
  &\leq c \sqrt{\e} \|\xi\|_X
  \| \eta \|_{L^2 \times L^2}
\end{split}
\]
and hence
\be
\la{eq:inq3}
  \|Q_\s \J L_\s \xi \|_{L^2 \times L^2}
  \leq c \sqrt{\e} \|\xi\|_X.
\end{equation}
Collecting~(\ref{eq:inq1}),~(\ref{eq:inq2}),
and~(\ref{eq:inq3}), we obtain a bound
on the right hand side of the effective dynamics law~(\ref{eq:ml}):
\be
\la{eq:eff}
  \| \p_t w_\s - Q_\s \J \cH'(w_\s) \|_{L^2 \times H^{-s}}
  \leq c(\sqrt{\e} + \|\xi\|_X + \| \dot \s \|_{Y_{1,0}})
  \|\xi\|_X.
\end{equation}

Finally, we translate~(\ref{eq:eff}) into parametric form,
in order to remove $\dot \s$ from the r.h.s., and
to see that it yields~(\ref{eq:mhlaw})
in the leading order.
For $\s \in C^1(\R;\Si_\e)$, we recall
\be
\la{eq:star}
  \p_t w_\s = \G_\s\dot \s
\end{equation}
where $\dot \s = (\dot z, \dot p, \p_t^z \chi, \p_t \z)$. Next,
using~(\ref{eq:exprQ}), we find \be \la{eq:sharp}
  Q_\s \J \cH'(w_\s) = \G_\s V_\s^{-1} \Lam_\s \cH'(w_\s).
\end{equation}
Now the definition of $\Lam_\s$,~(\ref{eq:lamdef}),
implies that
\be
\la{eq:sharp2}
  \Lam_\s \cH'(w_\s) = D_{\s}\cH(w_\s) = \p_\s h(\s).
\end{equation}
The last two equations yield
\be
\la{eq:star2}
  Q_\s \J \cH'(w_\s) = \G_\s V_\s^{-1} D_{\s}\cH(w_\s).
\end{equation}
Comparing~(\ref{eq:star}) with~(\ref{eq:star2}), we obtain
\[
  \p_t w_\s - Q_\s \J \cH'(w_\s)= \G_\s(\dot \s - V_\s^{-1} \p_\s
  h(\s)).
\]
Now the mapping properties~(\ref{eq:gamap}) and~(\ref{eq:lamap}),
% ?? with appropriate $r$ and $s$ here ??
and the fact that
$\Lam_\s \G_\s : Y_{1,1-s} \to Y_{-1,-s-1}$
is invertible,
% ?? why ??
imply that
\[
  \| \dot \s - V_\s^{-1} \p_\s h(\s) \|_{Y_{1,1-s}}
  \leq c \|\p_t w_\s - Q_\s \J \cH'(w_\s)\|_{L^2 \times H^{-s}}.
\]
Since $Dh(\s) = \lan \cH'(w_\s),D_\s w_\s \ran = O(\e^{1/2})$
(this estimate is part of Lemma~\ref{lem:approx} below),
(\ref{eq:eff}) implies that
\[
  \| \dot \s - V_\s^{-1} D_\s h(\s) \|_{Y_{1,1-s}}
  \leq c(\sqrt{\e} + \|\xi\|_X)\|\xi\|_X.
\]

Now using~(\ref{eq:dexpansion}) and~(\ref{eq:Vs}), we arrive at
\[
  V_\s^{-1} D_\s h(\s) =
  (\pb, \z, -\underline{\gamma}^{-1} \cdot \p_{\zb} W, 0)
  + O(\e\log^{1/2}(1/\e) (|\pb| + \|\zeta \|_2)),
\]
where we have used the notation
$\underline{\g}^{-1} \cdot \p_{\zb} W
= ( \g_{n_1}^{-1} \p_{z_1} W, \ldots, \g_{n_m}^{-1} \p_{z_m} W )$.
Combine this relation with the estimate above to obtain finally
\be
\la{eq:eff2}
\begin{split}
  \sum_{j=1}^m |\dot{z}_j(t) - p_j(t)| &+
  \sum_{j=1}^m |\dot{p}_j(t) + \g^{-1}_{n_j} \nabla_{z_j} W(\zb(t))|
  + \|\p_t^{\zb} \chi(t) - \z(t) \|_{H^1}
  + \|\p_t \z(t)\|_{H^{1-s}} \\
  &\leq c[(\sqrt{\e} + \|\xi\|_X)\|\xi\|_X + \e^{3/2}\log^{1/2}(1/\e)].
\end{split}
\end{equation}

%-------------------------------------------------------
\subsubsection{Energy estimates}

Our remaining task is to
control the remainder $\xi(t)$ for long times. The idea is similar
to techniques used to prove orbital stability of solitary waves in
Hamiltonian systems (see, eg, \cite{w,gss}): exploit conservation
of energy -- in this case {\it both} for the
PDE~(\ref{eq:hm}) {\it and} for the leading order effective
dynamics~(\ref{eq:mhlaw}) -- in order to control the fluctuations.
We begin with a Taylor expansion of the Hamiltonian: \be
\la{eq:taylor}
  \cH(w_\s + \xi) = \cH(w_\s) + \lan \cH'(w_\s), \xi \ran
  + \frac{1}{2} \lan \xi, L_\s \xi \ran
  + R_\s(\xi)
\end{equation}
(this equation defines $R_\s(\xi)$).
The following lemma, proved in Section~\ref{sec:coer},
allows us to control $\xi$ by the Hamiltonian.
\begin{lem}[coercivity]
\la{lem:coer}
For $\e$ sufficiently small,
$\s \in \Si_\e$, and $\xi \in ker Q_\s$,
\[
  \frac{1}{c} \|\xi\|_X \leq \lan \xi, L_\s \xi \ran
  \leq c \| \xi \|_X^2.
\]
\end{lem}

Using this lemma, together with conservation
of the Hamiltonian, in~(\ref{eq:taylor}), we obtain
\be
\la{eq:control}
\begin{split}
  \| \xi \|_X^2 &\leq c[\cH(w(0)) - \cH(w_\s) -
  \lan \cH'(w_\s), \xi \ran - R_\s(\xi) ] \\
  &= c[\cH(w_{\s(0)}) - \cH(w_\s)
  + \lan \cH'(w_{\s(0)}), \xi(0) \ran
  - \lan \cH'(w_\s), \xi \ran \\
  &\quad+ \frac{1}{2} \lan \xi(0), L_{\s(0)} \xi(0) \ran
  + R_{\s(0)}(\xi(0)) - R_\s(\xi)].
\end{split}
\end{equation}
The following lemma bounds the super-quadratic terms
on the right-hand side of~(\ref{eq:control}).
\begin{lem}[nonlinear estimate 2]
\la{lem:nl2}
For $\s \in \Si$,
\[
  |R_\s(\xi)| \leq c( \|\xi\|_X^3 + \|\xi\|_X^4 ).
\]
\end{lem}
This is proved is Section~\ref{sec:tay}.
To control the terms linear in $\xi$, we need another
key lemma:
\begin{lem}[approximate solution properties]
\la{lem:approx}
For $\s \in \Si_\e$, we have
\begin{enumerate}
\item
$\| \cH'(w_\s) \|_{H^1 \times L^2} \leq c \sqrt{\e}$
\item
$\|[\cH'(w_\s)]_1\|_{H^1} = \|\E_{GL}'(\vz)\|_{H^1}
\leq c \e \log^{1/4}(1/\e)$
\item
$\| \bQ_\s \J \cH'(w_\s) \|_{L^2 \times L^2}
\leq c \e \log^{1/4}(1/\e)$,
where $\bQ_\s := {\bf 1} - Q_\s$.
\end{enumerate}
\end{lem}
This lemma is proved in Section~\ref{sec:approx}. Using the third
statement of the lemma and~(\ref{eq:Q*}), we find
\be
\la{eq:bound}
\begin{split}
  |\lan \cH'(w_\s), \xi \ran|& = |\lan \cH'(w_\s), \bQ_\s \xi
  \ran|\\
  &= |\lan \J \bQ_\s \J \cH'(w_\s), \xi \ran|\\
  &\leq \| \bQ_\s \J \cH'(w_\s) \|_{L^2 \times L^2}
  \| \xi \|_X \\
  &\leq c \e \log^{1/4}(1/\e) \| \xi \|_X.
\end{split}
\end{equation}

Collecting Lemma~\ref{lem:coer}, Lemma~\ref{lem:nl2},
and~(\ref{eq:bound}) in~(\ref{eq:control})
(and remembering the intermediate assumption $\|\xi\|_X < 1$),
we obtain
\be
\la{eq:control2}
\begin{split}
  \|\xi\|_X^2 \leq c[
  &\cH(w_{\s(0)}) - \cH(w_\s) +
  (\e \log^{1/4}(1/\e) + \|\xi\|_X^2)\|\xi\|_X \\
  &+(\e \log^{1/4}(1/\e) + \|\xi(0)\|_X)\|\xi(0)\|_X ].
\end{split}
\end{equation}

%-----------------------------------------------------
\subsubsection{Approximate conservation of the reduced energy,
$\cH(w_\s)$}

It remains to control $\cH(w_{\s(0)}) - \cH(w_\s)$.
The estimate below involves a delicate estimate of the
contribution of the nonlinear terms.
\begin{prop}
Let $M(t) := \sup_{0 \leq s \leq t} \|\xi(s)\|_X$.  Then
\be
\la{eq:hamdiff}
  |\cH(v_{\s(0)}) - \cH(v_{\s(t)})| \leq
  ct \sqrt{\e} M(t)(\e\log^{1/2}(1/\e) + M(t)) + \sqrt{\e}M^2(t).
\end{equation}
\end{prop}
{\bf Proof}:
First, we differentiate in time and use the
effective dynamics law~(\ref{eq:ml}):
\be
\la{eq:enest}
\begin{split}
  \frac{d}{dt} \cH(w_\s)
  &= \lan \cH'(w_\s), \p_t w_\s \ran \\
  &= \lan \cH'(w_\s), Q_\s \J \cH'(w_\s)
  + Q_\s[\J L_\s \xi - \p_t \xi]
  + Q_\s \J N_\s(\xi) \ran \\
  &= \lan \cH'(w_\s), Q_\s[\J L_\s \xi - \p_t \xi] \ran
   + \lan \cH'(w_\s), Q_\s \J N_\s(\xi) \ran
\end{split}
\end{equation}
where we have used the fact that $(Q_\s \J)^* = -Q_\s \J$.

We start by estimating the first inner-product
on the right-hand side.
First we exploit the fact -- Lemma~\ref{lem:approx} part 2 --
that the first component of $\cH'(w_\s)$ is smaller than
the second: using
$Q_\s \p_t \xi = [Q_\s,\p_t]\xi$,
(\ref{eq:Q'}) and~(\ref{eq:inq3}), we have
\[
  |\lan [\cH'(w_\s)]_1,
  [Q_\s(\J L_\s \xi - \p_t\xi)]_1 \ran|
  \leq c \e \log^{1/4}(1/\e)
  \| \xi \|_X (\sqrt{\e}+\|\dot{\s}\|_{Y_{1,0}}).
\]
Combined with~(\ref{eq:eff2}), this yields
\[
  |\lan [\cH'(w_\s)]_1,
  [Q_\s(\J L_\s \xi - \p_t\xi)]_1 \ran|
  \leq c \e^{3/2} \log^{1/4}(1/\e) \| \xi \|_X.
\]

To deal with the second component, we have to
exploit a key cancellation.  This is expressed
in the following lemma, which can be considered
a refinement of both~(\ref{eq:inq2}) and
of~(\ref{eq:inq3}).
\begin{lem}
\la{lem:prec}
For $\s \in C^1(\R;\Si_\e)$ and $Q_\s \xi \equiv 0$,
we have
\[
\begin{split}
  |\lan [\cH'(w_\s)]_2, &
  [Q_\s(\J L_\s \xi - \p_t\xi)]_2 \ran|
  \leq c \sqrt{\e} [\e \log^{1/2}(1/\e) + |\dot \pb + \nabla_z W(\zb)| \\
  &+ \| \p_t \z \|_2 + \sqrt{\e}(|\dot{\zb}-\pb|
  + \|\p_t^{\zb} \chi - \z \|_{H^1})]\|\xi\|_X.
\end{split}
\]
\end{lem}
This lemma is proved in Section~\ref{sec:proj}.
Combining the above estimates
with~(\ref{eq:eff2}) yields
\be
\la{eq:dagger}
  |\lan \cH'(w_\s), Q_\s[\J L_\s \xi - \p_t \xi] \ran |
  \leq c [\e^{3/2} \log^{1/2}(1/\e) + \e \|\xi\|_X]\|\xi\|_X.
\end{equation}

Finally, we must control the second inner-product
on the right hand side of~(\ref{eq:enest}).
This is problematic since, so far, we have control
over (the second component of) the first factor only in $L^2$
(see Lemma~\ref{lem:approx}),
and over the second factor only in $H^{-s}$
for $s>0$ (see Lemma~\ref{lem:nl1}).
The solution is to isolate the worst term and
to use the detailed structure
of the equations to deal with it.

First we claim that $\cH'(w_\s)$ is of the form
\be
\la{eq:h'1}
  \cH'(w_\s) = \J H_\z^\s + \cH'_{rest}
\end{equation}
with
$H_\z^\s := \lan \z, \p_\chi w_\s \ran$,
and $\cH'_{rest}$ satisfying the estimate
\be
\la{eq:h'2}
  \| \cH'_{rest} \|_{H^{1-s} \times H^1}
  \leq c \sqrt{\e}.
\end{equation}
Indeed,~(\ref{eq:h'1})-~(\ref{eq:h'2}) is easily
obtained from the explicit expression
\be
\la{eq:h'exp}
  \cH'(w_\s) = (\E_{GL}'(\vz), \phi_\s),
\end{equation}
where $\phi_\s = p_{jk} \Tz_{jk} + \Gz_\z,$, and estimates
$\| \E_{GL}'(\vz) \|_{H^1} \leq c\e \log^{1/4}(1/\e)$
(Lemma~\ref{lem:approx}, part 2)
and $|p| + \|\z\|_{H^1} \leq \sqrt{\e}$.
Thus
\[
\begin{split}
  \lan \cH'(w_\s),Q_\s \J N_\s(\xi) \ran
  &= -\lan Q_\s [-H^\s_\z + \J \cH'_{rest}], N_\s(\xi) \ran \\
  & = \lan \Gz_\z, (N_\s)_1(\xi_1) \ran
  - \lan Q_\s \J \cH'_{rest}, N_\s(\xi) \ran.
\end{split}
\]
So by~(\ref{eq:h'2}), (\ref{eq:inq1}), and~(\ref{eq:Q*}),
we have
\be
\la{eq:dH1}
  |\lan \cH'(w_\s),Q_\s \J N_\s(\xi) \ran -
  \lan \Gz_\z, (N_\s)_1(\xi_1) \ran|
  \leq c \sqrt{\e} \|\xi\|_X^2.
\end{equation}

Next, we single out the worst term in the nonlinearity:
\[
  (N_\s)_1(\xi_1) = ( 0, Im(\bar{\xi}_{\psi} \nabla_{\Az} \xi_{\psi}))
  + N_{rest}
\]
where
\[
  \| N_{rest} \|_{H^{-s} \times L^2} \leq c \|\xi\|_X^2.
\]
Recall here that we are writing
$\xi = (\xi_1,\xi_2)$,
$\xi_1 = (\xi_\psi, \xi_A)$, and
$\xi_2 = (\xi_\pi, \xi_E)$.
Hence
\[
  | \lan G_\z, N_{rest} \ran | \leq c \|\z\|_{H^1} \|\xi\|_X^2.
\]
Then in light of~(\ref{eq:dH1}), we have
\be
\la{eq:dH2}
  |\lan \cH'(w_\s),Q_\s \J N_\s(\xi) \ran -
  \lan \nabla \z,
  Im(\bar{\xi}_{\psi} \nabla_{\Az} \xi_\psi) \ran|
  \leq c \sqrt{\e} \|\xi\|_X^2.
\end{equation}

It remains to estimate $\lan \nabla \cdot \z, Im(\bar{\xi}_{\psi}
\nabla_{\Az} \xi_\psi) \ran$ (note that the estimates available to
us so far control the first factor in $L^2$ (and no better) and
just fail to control the second in $L^2$). The key is to recognize
this quantity as (essentially) a time derivative.  Using the basic
equation~(\ref{eq:exp}) and~(\ref{eq:h'exp})
compute
\[
\begin{split}
  \frac{d}{dt} \lan \z, Im(\bar{\xi}_{\psi} \xi_\pi) \ran
  &= \lan \p_t \z, Im(\bar{\xi}_{\psi} \xi_\pi) \ran
  + \lan \z, Im( \bar{\xi}_{\psi}[-\p_t [\phi_\s]_\pi
  + [\E'_{GL}(\vz)]_\psi \\
  &+ [\E''_{GL}(\vz) \xi_1]_\psi + [(N_\s)_1]_\psi ]
  + \bar{\xi}_\pi[\p_t \psiz - [\phi_\s]_\psi]) \ran.
\end{split}
\]
We estimate each term on the RHS as follows:
\[
  | \lan \p_t \z, Im(\bar{\xi}_{\psi} \xi_\pi) \ran |
  \leq c \| \p_t \z \|_{H^{1-s}} \|\xi\|_X^2
\]
% ?? true ??
\[
  | \lan \z, Im(\bar{\xi}_{\psi} [(N_\s)_1]_\psi) \ran |
  \leq c \|\z\|_{H^1}(\|\xi\|_X^3+\|\xi\|_X^4)
\]
\[
  | \lan \z, Im(\bar{\xi}_{\psi} [\E_{GL}'(\vz)]_\psi) \ran |
  \leq c \e \log^{1/4}(1/\e) \|\z\|_{H^1} \|\xi\|_X
\]
\[
\begin{split}
  | \lan \z, Im(\bar{\xi}_{\psi} \p_t [\phi_\s]_\pi) \ran |
  &\leq c (|\dot{\pb}| + (|\pb|+\|\p_t \z\|_2)
  (|\dot{\zb}|+\|\p_t^{\zb} \chi\|_2) \\
  &\quad + \| \p_t \z \|_2) \|\z\|_{H^1}\|\xi\|_X
\end{split}
\]
and
\[
\begin{split}
  | \lan \z, & Im(\bar{\xi}_\pi[\p_t \psiz - [\phi_\s]_\psi]) \ran | \\
  & = | \lan \z, Im(\bar{\xi}_\pi
  [(\pb-\dot{\zb})_{jk} [\Tz_{jk}]_\psi +
  [\Gz_{\p_t^{\zb} \chi - \z}]_\psi]) \ran | \\
  &\leq c(|\pb-\dot{\zb}| + \| \p_t^{\zb} \chi - \z \|_{H^1})
  \|\z\|_{H^1} \|\xi\|_X.
\end{split}
\]
We write
\[
  [\E_{GL}''(\vz) \xi_1]_{\psi} =
  -\Delta_{\Az} \xi_\psi + \E''_{rest}
\]
with
\[
  | \lan \z, Im(\bar{\xi}_{\psi} \E''_{rest}) \ran |
  \leq c \|\z\|_{H^1} \|\xi\|_X^2.
\]

Collecting these estimates yields
\be
\la{eq:collect}
\begin{split}
  |\frac{d}{dt} &\lan \z, Im(\bar{\xi}_{\psi} \xi_\pi) \ran
  + \lan \z, Im(\bar{\xi}_{\psi} \Delta_{\Az} \xi_\psi) \ran | \\
  &\leq c \|\xi\|_X( \| \p_t \z \|_{H^{1-s}} \|\xi\|_X
  + \sqrt{\e}[\|\xi\|_X + \e\log^{1/4}(1/\e) + |\dot{\pb}| \\
  &\quad + \sqrt{\e}(|\dot{\zb}|+\|\p_t^{\zb} \chi\|_2)
  + |\pb - \dot{\zb}| + \|\p_t^{\zb} \chi - \z\|_{H^1}]) \\
  &\leq c \sqrt{\e}\|\xi\|_X (\|\xi\|_X + \e \log^{1/4}(1/\e)),
\end{split}
\end{equation}
using~(\ref{eq:eff2}).  Noting that
\[
  -\lan \z, Im(\bar{\xi}_{\psi} \Delta_{\Az} \xi_\psi) \ran
  = \lan \nabla \z, Im(\bar{\xi}_{\psi}
  \nabla_{\Az} \xi_\psi) \ran,
\]
and combining~(\ref{eq:enest}),~(\ref{eq:dagger}),~(\ref{eq:dH2}),
and~(\ref{eq:collect}), we have
\[
  |\frac{d}{dt} [\cH(w_\s) + \lan \z,
  Im(\bar{\xi}_{\psi} \xi_\pi) \ran]|
  \leq c\sqrt{\e}\|\xi\|_X(\|\xi\|_X+\e\log^{1/4}(1/\e)).
\]
Integrating this in time, and defining
$M(t) := \sup_{0 \leq s \leq t} \|\xi(s)\|_X$,
leads to~(\ref{eq:hamdiff}).
$\Box$

Returning to~(\ref{eq:control2}), we obtain
\[
\begin{split}
  \|\xi(t)\|_X^2 &\leq c[
  t\sqrt{\e}M(t)(M(t)+\e\log^{1/4}(1/\e))
  + M(t)(\e\log^{1/2}(1/\e)+M^2(t)) \\
  &\quad + \sqrt{\e} M^2(t) M(0)(\e\log^{1/4}(1/\e)+M(0))].
\end{split}
\]
It now follows that there is a constant $\tau' > 0$, such that for
$0 \leq t \leq \min \left( \frac{\tau'}{\sqrt{\e}}, T_1 \right)$
(recall $T_1$ is the time of first exit of $w(t)$
from $U_\del$), we have
\be
\la{eq:xibound}
  \|\xi(t)\|_X < c(\e\log^{1/2}(1/\e) + \|\xi(0)\|_X).
\end{equation}
In particular, the intermediate assumption~(\ref{eq:xismall})
is justified.

%------------------------------------------------------
\subsubsection{A priori momentum bound}

We wish to iterate the above argument to extend the time interval.
The problem is that the vortex velocities can, in principle, grow
to size $\e t >> \sqrt{\e}$ if $t >> 1/\sqrt{\e}$ (this would mean
we leave the manifold $M_{mv}$, and many of the above estimates fail).
We show here that this does not happen. To this end we use the
approximate conservation of the reduced (vortex) energy $H(w_\s)$,
together with the repulsivity of the interaction energy.
Indeed, since we are in the ``repulsive'' case
($\lambda > 1/2$ and $n_1 = \cdots = n_m = \pm 1$),
we have the following lemma, which is proved
in Section~\ref{sec:reden}.
\begin{lem}[interaction energy]
\la{lem:repulse}
For $R(\zb)$ large,
\be
\la{eq:repulse}
  W(\zb) = \sum_{j \not= k} n_j n_k c_{jk}
  \frac{e^{-|z_j-z_k|}}{\sqrt{|z_j-z_k|}}
  + o(e^{-R(\zb)}/\sqrt{R(\zb)}).
\end{equation}
Here $c_{jk} > 0$ are constants.
\end{lem}
\begin{rem}
One can see from this expression that like-signed
vortices repel, while opposite-signed vortices
attract.
\end{rem}

Conservation of energy for the PDE,
$\cH(w(t)) = \cH(w(0))$, together with
the decomposition~(\ref{eq:splitter})
and a Taylor expansion yields
\[
\begin{split}
  \cH(w_\s) - \cH(w_{\s(0)})
  &= \lan \cH'(w_\s), \xi \ran
  - \lan \cH'(w_{\s(0)}), \xi(0) \ran \\
  & \quad + O(\|\xi\|_X^2 + \|\xi(0)\|_X^2).
\end{split}
\]
We saw above (in~(\ref{eq:bound})) that
$|\lan \cH'(w_\s), \xi \ran| \leq c \e \log^{1/4}(1/\e) \|\xi\|_X$,
which gives
\[
  \cH(w_\s)-\cH(w_{\s(0)})
  \leq c((\|\xi\|_X+\|\xi(0)\|_X)
  (\e \log^{1/4}(1/\e)+\|\xi\|_X + \|\xi(0)\|_X)).
\]
By estimates~(\ref{eq:expansion}) and $K\ge c$ for some $c>0$, we
have, for $\e$ sufficiently small,
\[
  W(\zb) + |\pb|^2 + \|\z\|_{H^1}^2
  < \alpha[\cH(w_\s) - \sum_{j=1}^m E^{(n_j)}]
\]
for some constant $\alpha > 0$.  In light of Lemma~\ref{lem:repulse},
the assumptions on the initial conditions in Theorem~\ref{thm:mh}
imply $\cH(w_{\s(0)}) - \sum_{j=1}^m E^{(n_j)} < c' d_0 \e$,
for some $c'$. We choose $d_0 < \frac{\alpha'}{2c'\alpha}$,
where $\alpha'$ is a constant to be chosen below.
So provided
\be
\la{eq:condition}
  c((\|\xi(t)\|_X + \|\xi(0)\|_X)
  (\e \log^{1/4}(1/\e) + \|\xi(t)\|_X + \|\xi(0)\|_X))
  < \frac{\alpha'}{2\alpha} \e,
\end{equation}
we have
$\alpha[\cH(w_\s) - \sum_{j=1}^m E^{(n_j)}] < \alpha' \e$, and therefore
\be
\la{eq:enbound}
  W(\zb) + |\pb|^2 + \| \z \|_{H^1}^2 < \alpha' \e.
\end{equation}

So by~(\ref{eq:repulse}), if $\al' < \min(c_{jk})$, then
as long as condition~(\ref{eq:condition}) holds,
$|p|^2 + \| \z \|_{H^1}^2 < \e$, and
$R(\zb)e^{-R(\zb)} < \e$.  Hence $\s \in \Si_\e$,
and $w_\s \in M_{mv}$.

In particular, this estimate shows that
$T_1 > \tau'/\sqrt{\e}$.
Hence, we have shown:
\begin{lem}
\la{lem:iter}
There are $\tau' > 0$ and $d > 0$, such that
inequality~(\ref{eq:xibound}) holds
for $0 \leq t \leq \tau'/\sqrt{\e}$,
provided
\be
\la{eq:iter}
  \|\xi(0)\|_X + \|\xi(t)\|_X < d \sqrt{\e}.
\end{equation}
\end{lem}

%------------------------------------------------------
\subsubsection{Iteration}

We may iterate Lemma~\ref{lem:iter} for as long as the conditions
$\s \in \Si_\e$, and $\|\xi(t)\|_X < d \sqrt{\e}$ hold. Iterating
$N$ times starting with $\xi(0)$ and satisfying $\|\xi(0)\|_X\ \le
d_0\sqrt{\e}$ yields
\[
  \| \xi(t) \|_X \leq C c^N \e \log^{1/2}(1/\e)
  \quad \mbox{ for } \quad
  0 \leq t \leq \tau' N/\sqrt{\e},
\]
where $C$ is another constant.
The condition~(\ref{eq:iter}) limiting the number of iterations,
ensures both that~(\ref{eq:condition}) holds
(so that $\s \in \Si_\e$ remains true), and that
the remainder in the effective dynamics law is sub-leading order.
Thus we can take $c^N = \al(\e)/\sqrt{\e}$ for any
$\al_0 \sqrt{\e} < \al(\e) << \log^{-1/2}(1/\e)$ (with $\al_0 > 1$).
This gives a total time interval of length
\[
  T = \frac{\ta}{\sqrt{\e}}
  \log\left(\frac{\al(\e)}{\sqrt{\e}}\right),
\]
where $\ta = \ta'/\log c$,
over which we have the bound
\be
\la{eq:xibound2}
  \| \xi(t) \|_X < C \al(\e) \sqrt{\e} \log^{1/2}(1/\e)
\end{equation}
for $0 \leq t \leq T$.

Finally, equation~(\ref{eq:eff2}) implies
\[
\begin{split}
  |\dot{\zb}(t) - \pb(t)| &+
  |\dot{\pb}(t) + \underline{\gamma}^{-1} \cdot \nabla_{\zb} W(\zb)| +
  \|\p_t^{\zb} \chi(t) - \z(t) \|_{H^{1-s}} +
  \|\p_t \z(t) \|_{H^1} \\
  &\leq c \sqrt{\e} \|\xi\|_X
  \leq c \e \al(\e) \log^{1/2}(1/\e) = o(\e).
\end{split}
\]
This completes the proof of Theorem~\ref{thm:mh}.
$\Box$

%----------------------------------------------------------------
\subsection{Effective dynamics of vortices: superconductor model}
\label{sec:gf}

Here we just sketch the proof of  Theorem~\ref{thm:gf}
since it proceeds as above. The important difference is that
we can control the remainder for all times.

The set-up is as follows. For the superconductor model our
manifold of multi-vortex configurations is taken to be
\[
  M_{mv} = \{ \vz \; | \; e^{-R(\zb)}/\sqrt{R(\zb)} < \e, \;
  \chi \in H^2_{\zb}(\R^2;\R) \}.
\]
A solution $u(t) = (\psi(t), A(t))$ of~(\ref{eq:gf}) is decomposed
as
\[
  u(t) = \vzt + \xi(t)
\]
with $\Pz \xi \equiv 0$, where
$\Pz$ denotes the orthogonal projection
from $H^1$ onto the tangent space $T_{\vz} M_{mv}$.
Substituting this into~(\ref{eq:gf}) yields
\[
  \p_t \vz + \p_t \xi = -[\E_{GL}'(\vz)
  + \Lz \xi - N_{\vz}(\xi)]
\]
where $\Lz := \E_{GL}''(\vz)$.
The equation governing the effective dynamics of
$\zb(t)$ and $\chi(t)$ is derived by applying the
projection $\Pz$ to this:
\[
  \p_t \vz + \Pz \E_{GL}'(\vz)
  = -\Pz[\Lz \xi - N_{\vz}(\xi) + \p_t \xi].
\]
To estimate the RHS, we have the following properties
\[
  \| \Pz \Lz \xi \|_{L^2}
  \leq c \e \log^{1/2}(1/\e) \|\xi\|_{H^1}
\]
\[
  \| \Pz N_{\vz}(\xi) \|_{H^{-s}}
  \leq c(\|\xi\|_{H^1}^2+\|\xi\|_{H^1}^3)
\]
\[
  \| \Pz \p_t \xi \|_{H^1}
  \leq c(|\dot{\zb}| + \|\p_t^{\zb} \chi\|_{L^2})\|\xi\|_{H^1}.
\]
Combining these yields
\[
  \| \p_t \vz + \Pz \E_{GL}'(\vz) \|_{H^{-s}}
  \leq c(\e\log^{1/2}(1/\e) + |\dot{\zb}| + \|\p_t^{\zb} \chi\|_2 +
  \|\xi\|_{H^1})\|\xi\|_{H^1}
\]
which in parametric form reads
(recall the notation
$\g \dot{\zb} := (\g_{n_1} \dot{z}_1, \ldots, \g_{n_m} \dot{z}_m)$)
\be
\la{eq:effgf}
  |\g \dot{\zb} + \nabla_z W(\zb)| + \|\p_t^{\zb} \chi \|_{H^{1-s}}
  \leq c(\e\log^{1/2}(1/\e) + \|\xi\|_{H^1})\|\xi\|_{H^1}.
\end{equation}

In order to control $\|\xi\|_{H^1}$ for {\it all} time we need
\begin{lem}
\label{lem:co2}
There is $\del > 0$ such that for $R(\zb)$ sufficiently large,
\[
  \lan \Lz \xi, \Lz \xi \ran > \del \|\xi\|_{H^2}^2.
\]
\end{lem}
This lemma is proved in Section~\ref{sec:coer}.

We use the fact that the main part of the energy difference,
$\E(\vz+\xi)-\E(\vz)$, namely
$\frac{1}{2} \lan \xi, \Lz \xi \ran$, is
a decaying quantity.  Compute
\[
\begin{split}
  \frac{d}{dt} \lan \xi, \Lz \xi \ran
  &= 2 \lan \p_t \xi, \Lz \xi \ran
    + \lan [\p_t, \Lz] \xi, \xi \ran \\
  &= \lan -\E'(\vz) - \Lz\xi - N_{\vz}(\xi), \Lz \xi \ran \\
  & \quad + O((|\dot{\zb}|+\|\p_t \chi \|_2)\|\xi\|_{H^1}^2).
\end{split}
\]
Now we use
\[
  |\lan \E'(\vz), \Lz \xi \ran| < c \e \log^{1/4}(1/\e)
  \|\xi\|_{H^1},
\]
\[
  |\lan  N_{\vz}(\xi), \Lz \xi \ran| <
  c \|\xi\|_{H^2}^2(\|\xi\|_{H^1}+\|\xi\|_{H^1}^2),
\]
and
\[
  \lan \xi, \Lz \xi \ran \leq \beta \|\xi\|_{H^1}^2,
\]
together with Lemma~\ref{lem:co2}, to obtain
\[
\begin{split}
  (\frac{d}{dt} + \frac{\del}{2\beta})
  \lan \xi, \Lz \xi \ran &\leq
  \|\xi\|_{H^2}^2
  [ c(\|\xi\|_{H^1}+\|\xi\|_{H^1}^2  \\
  & \quad +|\dot{\zb}|+\|\p_t \chi\|_2)
  - \del/2 ] + c \e \log^{1/4}(1/\e)\| \xi \|_{H^1}.
\end{split}
\]
So as long as
\be
\la{eq:temp}
  c(\|\xi\|_{H^1}+\|\xi\|_{H^1}^2
  + |\dot{\zb}|+\|\p_t \chi \|_2) \leq \del/2,
\end{equation}
we have
\[
  \frac{d}{dt} (e^{(\del/2\beta)t} \lan \xi, \Lz \xi \ran)
  \leq c \e \log^{1/4}(1/\e) e^{(\del/\beta)t} \|\xi\|_{H^1}.
\]
Setting $M(t) := \sup_{0 \leq s \leq t} \|\xi(s)\|_{H^1}$
and integrating in time leads to
\[
  \lan \xi, \Lz \xi \ran \leq
  e^{-(\del/2\beta)t} \lan \xi(0), \Lzo \xi(0) \ran
  + c \e \log^{1/4}(1/\e) M(t).
\]
Finally using
\[
  c \lan \xi, \Lz \xi \ran <
  \|\xi\|_{H^1}^2 <
  \frac{1}{\g} \lan \xi, \Lz \xi \ran,
\]
we find
\[
  M^2(t) \leq c[e^{-(\del/2\beta)t} M^2(0) +
  \e \log^{1/4}(1/\e) M(t)]
\]
and so if $M(0) = O(\e\log^{1/4}(1/\e))$ we have
\[
  \| \xi(t) \|_{H^1} < c \e \log^{1/4}(1/\e)
\]
for all $t$, as long as~(\ref{eq:temp})
and $e^{-R(\zb)}/\sqrt{R(\zb)} < \e$ hold.
By~(\ref{eq:effgf}), we see
\[
  |\dot{\zb} + \nabla_{\zb} W(\zb)| +
  \| \p_t^{\zb} \chi \|_{H^{1-s}} \leq c \e^2 \log^{3/4}(1/\e)
\]
so the intermediate assumption~(\ref{eq:temp})
is justified.  Finally, in the repulsive case,
$\E(u) - \sum_{j=1}^m E^{(n_j)} \leq c\e$ implies
$e^{-R(\e)}/\sqrt{R(\zb)} < \e$ holds for all $t$.
$\Box$

%%%%%%%%%%%%%%%%%%%%%%%%%%%%%%%%%%%%%%%%%%%%%%%%%
%%%%%%%%%%%%%%%%    Properties    %%%%%%%%%%%%%%%
%%%%%%%%%%%%%%%%%%%%%%%%%%%%%%%%%%%%%%%%%%%%%%%%%

\section{Key properties}
\la{sec:prop}

In this section we prove the lemmas used in the proofs of
Theorems~\ref{thm:gf} and~\ref{thm:mh}.

%--------------------------------------------------
\subsection{Approximate static solution property}
\la{sec:approx}

{\bf Proof of Lemma~\ref{lem:approx}}:
The main fact we use here is that since we consider
the Type-II regime ($\lam > 1/2$),
the effects of the magnetic field
and current dominate those of the order parameter
at large distances.

In what follows, a subindex $k$ will
denote an equivariant field component, of degree $n_k$,
centred at $z_k$:
eg, $\psi_k := \psi^{\nk}(\cdot-z_k)$,
$(\nabla_A \psi)_k = \nabla_{A^{\nk}(\cdot-z_k)} \psi^{\nk}(\cdot-z_k)$,
etc.

We first prove
\begin{lem}
\la{prop:de2}
\be
\la{eq:de2}
  \| \E_{GL}'(\vz) \|_2 \leq c e^{-R(\zb)}/R^{1/4}(\zb).
\end{equation}
\end{lem}
{\bf Proof:}
The proof is a computation using the fact that
$u^{\nj} = (\psi^{\nj},A^{\nj})$
satisfies the Ginzburg-Landau equations,
together with the exponential decay~(\ref{eq:decay}).
We start with
\[
  [\E_{GL}'(\vz)]_\psi = -\Delta_{\Az} \psiz + \lam(|\psiz|^2-1)\psiz.
\]
Using gauge covariance and the covariant product rule, we find
\be
\label{eq:11}
  \Delta_{\Az} \psiz = e^{i \chi} \left[
  \sum_j ( \prod_{k \not= j} \psi_k )
  (\Delta_A \psi)_j +
  \sum_{j \not= k} ( \prod_{l \not= j,k} \psi_l)
  (\dA \psi)_j \cdot (\dA \psi)_k \right].
\end{equation}
A little computation plus~(\ref{eq:decay}) yields
\be
\label{eq:12}
  | (\prod_{j=1}^m f_j^2) - 1 - \sum_{j=1}^m (f_j^2-1) |
  \leq c \sum_{j \not= k}
  e^{-m_{\lam}(|x-z_j| + |x-z_k|)}.
\end{equation}

Using~(\ref{eq:11}) and~(\ref{eq:12}),
together with the fact that $u^{\nj}$
solves the Ginzburg-Landau equations,
we arrive at
\[
  |[\E_{GL}'(\vz)]_\psi(x) - [E^{(\zb,\chi)}]_\psi(x)|
  \leq c \sum_{j \not= k} e^{-m_{\lam}(|x-z_j| + |x-z_k|)}
\]
where
\[
  E^{(\zb,\chi)}_\psi :=
  -e^{i \chi} \sum_{j \not= k}
  (\prod_{l \not= j,k} \psi_l)
  (\dA \psi)_j \cdot (\dA \psi)_k.
\]
Using Lemma~\ref{lem:*} (in Appendix 3, Section~\ref{sec:app3})
with $\al = \beta = 2m_\lam > 2$, $\gamma = \del = 0$, we obtain
\be \la{eq:piece1}
  \| [\E_{GL}'(\vz)]_\psi - [E^{(\zb,\chi)}]_\psi \|_2
  \leq c e^{-m_\lam R(\zb)} R(\zb)^{3/2}
  << e^{-R(\zb)}/\sqrt{R(\zb)}.
\end{equation}

We turn now to
\[
  [\E_{GL}'(\vz)]_A = curl \Bz - \jz.
\]
Observing that
$curl \Bz = \sum_{j=1}^m curl B_j$, and
\[
  \jz = \sum_{j=1}^m j_j + \sum_{j=1}^m
  (\prod_{k \not= j} f_k^2 - 1) j_j,
\]
using the Ginzburg-Landau equation
$curl B_j - j_j = 0$, invoking an equation
similar to~(\ref{eq:12}) for $\prod_{k \not= j} f_k^2 - 1$,
and using~(\ref{eq:decay}), we arrive at
\[
  |[\E_{GL}'(\vz)]_A(x) - E^{(\zb,\chi)}_A(x)|
  \leq c \sum_{j,k,l \mbox{ distinct}}
  e^{-m_{\lam}(|x-z_j| + |x-z_k|) - |x-z_l|}
\]
where
\[
  E^{(\zb,\chi)}_A := \sum_{j \not= k} (1-f_j^2) j_k.
\]
Estimating as above gives
\be
\la{eq:piece2}
  \| [\E_{GL}'(\vz)]_A - E^{(\zb,\chi)}_A \|_2
  \leq c e^{-m_\lam R(\zb)} R(\zb)^{3/2}
  << e^{-R(\zb)}/\sqrt{R(\zb)}.
\end{equation}

Using~(\ref{eq:decay}) again, we obtain
the following pointwise estimate for
$E^{(\zb,\chi)}=(E_\psi^{(\zb,\chi)},E_A^{(\zb,\chi)})$:
\[
  |E^{(\zb,\chi)}| \leq c \sum_{k \not= j}
  \frac{e^{-|x-z_j|}}{(1+|x-z_j|)^{1/2}}
  \frac{e^{-|x-z_k|}}{(1+|x-z_k|)^{1/2}}.
\]
Applying Lemma~\ref{lem:*} with $\alpha=\beta=2$
and $\gamma=\delta=1$ yields
\be
\la{eq:Ebound}
  \| E^{(\zb,\chi)} \|_2
  \leq c e^{-R(\zb)}/R(\zb)^{1/4}.
\end{equation}
Then~(\ref{eq:piece1})-(\ref{eq:Ebound}) yield~(\ref{eq:de2}).
$\Box$

Now we consider the
manifold of approximate solutions
for the Higgs model equations.
Parts 1 and 2 of Lemma~\ref{lem:approx} follow
immediately from the expression
(cf.~(\ref{eq:h'exp}))
\[
  \cH'(w_\s) = (\E_{GL}'(\vz),p_{jk} \Tz_{jk} + \Gz_\z),
\]
together with Lemma~\ref{prop:de2},
and the fact that $\s \in \Si_\e$ implies
$|\pb|+\|\z\|_{H^1} < \e$ and
$e^{-R(\zb)}/\sqrt{R(\zb)} < \e$,
which implies $e^{-R(\zb)}/R(\zb)^{1/4} < c \e \log^{1/4}(1/\e)$.
The refined statement, part 3 of Lemma~\ref{lem:approx},
follows from the fact that for $\s \in \Si_\e$,
\[
  Q_\s = \left( \ba{cc}
  \Pz & 0 \\ 0 & \Pz
  \ea \right)
  + O(\sqrt{\e})
\]
where $\Pz$ denotes the orthogonal projection
onto the span of $\Tz_{jk}$ and $\Gz_{\gamma}$
(see Eqn.~(\ref{eq:Qform}) of Appendix 2) and so,
since $\bPz[p_{jk} \Tz_{jk} + \Gz_\z] = 0$,
\[
  \| \bQ_\s \J \cH'(w_\s) \|_{L^2 \times L^2}
  \leq c\e \log^{1/4}(1/\e)
\]
as required.
This completes the proof of Lemma~\ref{lem:approx}.
$\Box$

%-------------------------------------------------------
\subsection{Inter-vortex interaction}
\la{sec:reden}

The reduced energy is a function of the
vortex positions alone:
\[
  W(\zb) := \E_{GL}(\vz) - \sum_{j=1}^m E^{\nj}.
\]
In this section, we compute -- to leading order
in the vortex separation --
$W(\zb)$, and $\nabla W(\zb)$, the inter-vortex force
entering the effective vortex dynamic laws.

\medskip

\noindent
{\bf Proof of Lemma~\ref{lem:repulse}}:
Noting that $B = curl A$, we re-write
the Ginzburg-Landau energy as
\[
  \E_{GL}(A,\psi) = \frac{1}{2} \int_{\R^2}
  \{ |\nabla_A \psi|^2 + B^2 + \frac{\lam}{2}(|\psi|^2-1)^2 \}.
\]
For $(\psi,A) = (\psiz,\Az)$, we have
\[
  \nabla_A \psi =
  e^{i\chi} \sum_{j=1}^m (\prod_{k \not= j} \psi_k)
  (\nabla_A \psi)_j,
\]
and $B = B_1 + \cdots + B_m$.  So plugging into $\E_{GL}$ and
using the notation $j_l := Im(\bar{\psi}_l \nabla_{A_l} \psi_l)$
and $f_l := |\psi_l|$, we find
\[
  \E_{GL}(\psiz,\Az)  = \sum_{j=1}^m E^{\nj}
  + LO + Rem
\]
where
\[
  LO := \frac{1}{2} \sum_{l \not= k} \int_{\R^2}
  [ j_l \cdot j_k + B_l B_k ]
\]
and
\[
\begin{split}
  Rem &=
  \frac{1}{2} \sum_{j=1}^m \int (\prod_{k \not= j} f_k^2-1)
  |(\nabla_A \psi)_j|^2
  + \frac{1}{2} \sum_{j \not= l} \int (\prod_{k \not= j,l} f_k^2)
  [Re(\bar{\psi}\nabla_A \psi)]_j [Re(\bar{\psi}\nabla_A \psi)]_l \\
  &\quad + \frac{1}{2} \sum_{j \not= l} \int (\prod_{k \not= j,l} f_k^2-1)
  j_k \cdot j_l \\
  &\quad + \frac{\lam}{4} \int [
  \sum_{j \not = l} (f_j^2-1)(f_l^2-1) +
  \sum_{j \not= l \not= k} (f_j^2-1)(f_l^2-1)(f_k^2-1) + \cdots].
\end{split}
\]
For each term in $Rem$, the integrand is bounded by
$e^{-(\min(m_{\lam},2)|x-z_j|+m_{\lam}|x-z_k|)}$ or
$e^{-(m_{\lam}|x-z_k|+|x-z_j|+|x-z_l|)}$, and so, after integration,
is $<< e^{-R(\zb)}/\sqrt{R(\zb)}$
(using Lemma~\ref{lem:*} and $m_{\lam} > 1$).
Using the Ginzburg-Landau equation
$curl B = j$, we can re-write the leading-order term as
\[
  LO = \frac{1}{2} \sum_{l \not= k} \int_{\R^2}
  [ B_l (-\Delta + 1)B_k ].
\]
A computation gives
$(-\Delta+1)B = n(2(1-a)ff' + a'(1-f^2))/r > 0$.
By~(\ref{eq:decay}),
$|(-\Delta+1)B(x)| < c e^{-m_{\lam} |x|}$,
and $B_n(x) = c_n n e^{-|x|}/\sqrt{|x|}[1+O(1/|x|)]$,
$c_n > 0$. Applying Lemma~\ref{lem:**} yields
\[
\begin{split}
  LO &= \frac{1}{2} \sum_{l \not=k} c_l n_l n_k
  \frac{e^{-|z_l-z_k|}}{\sqrt{|z_l-z_k|}}
  \int_{\R^2} e^{x \cdot (z_l-z_k)/|z_l-z_k|}
  (2(1-a_k)f_kf_k'+a_k'(1-f_k^2))/r dx \\
  & \quad + o\left(\frac{e^{-R(\zb)}}{\sqrt{R(\zb)}} \right).
\end{split}
\]
Lemma~\ref{lem:repulse} follows.
$\Box$.

Now we turn to the estimate of the force:
\begin{lem}
\label{prop:force}
We have
\be
\la{eq:forceasy}
  \nabla_{z_l} W(\zb) =  \sum_{j \not= l} n_j n_l
  C_{jl} \frac{e^{-|z_j-z_l|}}{\sqrt{|z_j-z_l|}}
  \frac{z_j-z_l}{|z_j-z_l|}
  + o(e^{-R(\zb)}/\sqrt{R(\zb)})
\end{equation}
as $R(\zb) \ra \infty$.
Here $C_{jl} > 0$ are constants.
\end{lem}
{\bf Proof:}
By the definition of $W(\zb)$,
$\nabla_{z_{lm}} W(\zb) = \lan \E_{GL}'(\vz), \Tz_{lm} \ran$.
Equations~(\ref{eq:piece1}) and~(\ref{eq:piece2}) imply that
\be
\la{eq:force}
  \nabla_{z_{lm}} W(\zb) =
  \lan E^{(\zb,\chi)}, \Tz_{lm} \ran
  + o(e^{-R(\zb)}/\sqrt{R(\zb)}),
\end{equation}
where $E^{(\zb,\chi)}$ is defined in  the proof of
Lemma~\ref{prop:de2}. We first compute
\[
  \lan E_{\psi}^{(\zb,\chi)}, [\Tz_{lm}]_{\psi} \ran =
  \sum_{j \not= k} \al_{lm}^{jk}
\]
where (recall the notation
$\psi_k(x) := \psi^{\nk}(x-z_k)$, etc.)
\[
  \al_{lm}^{jk} =
  \lan (\prod_{r \not= j,k} \psi_r) (\dA \psi)_j
  \cdot (\dA \psi)_k,
  (\prod_{t \not= l} \psi_t) ([\nabla_A]_m \psi)_l \ran.
\]
First, we note that $\al_{lm}^{jk} = \al_{lm}^{kj}$.
Second, we use~(\ref{eq:decay}) to conclude that
if $l \not= j$ and $l \not= k$, then
$| \al_{lm}^{jk} | << e^{-R(\zb)}/\sqrt{R(\zb)}$.
It remains to compute $\al_{jm}^{jk}$.  We rewrite to get
\[
  \al_{jm}^{jk} = \sum_s \int ( \prod_{r \not= j,k} f_r^2 )
  Re [ \overline{([\nabla_A]_m \psi)_j} ([\nabla_A]_s \psi)_j
  (\overline{\psi} [\nabla_A]_s \psi)_k ],
\]
and use
\[
  |\prod_{r \not= j,k} f_r^2 - 1| \leq c
  e^{-m_{\lam} \cdot \max(|x-z_j|,|x-z_k|)}
\]
to conclude that $\al_{jm}^{jk} = \alt_{jm}^{jk}
+ o(e^{-R(\zb)}/\sqrt{R(\zb)})$,
where
\[
  \alt_{jm}^{jk} = \sum_s Re \int \overline{([\nabla_A]_m \psi)_j}
  ([\nabla_A]_s \psi)_j
  (\overline{\psi} [\nabla_A]_s \psi)_k.
\]
Writing everything out in terms of the vortex profiles
$f_j$ and $a_j$ and taking the real part,
we find (applying Lemma~\ref{lem:*} again) that
\[
  \alt_{lm}^{lk} = - \sum_s \int \left[ \frac{n(1-a)}{r} (J\hat{x})_s \right]_k
  \left[ \frac{n(1-a)ff'}{r}
  [\hat{x}_m (J\hat{x})_s - (J\hat{x})_m\hat{x}_s] \right]_j
  + o(e^{-R(\zb)}/\sqrt{R(\zb)}).
\]
Now using the fact that
$[\hat{x}_m (J\hat{x})_s - (J\hat{x})_m\hat{x}_s]$ equals
$0$ if $s=m$, $-1$ if $(s,m)=(1,2)$, and $1$ if $(s,m)=(2,1)$,
and summing over $s$, we arrive at
\[
  \alt_{lm}^{lk} = - \int \left[ \frac{n(1-a)}{r} \right]_k
  \left[ \frac{n(1-a)ff'}{r} \right]_j
  (\widehat{x-z_k})_m
  + o(e^{-R(\zb)}/\sqrt{R(\zb)}).
\]
Now we apply (a slight variant of) Lemma~\ref{lem:**} to obtain
\[
\begin{split}
  \alt_{jm}^{jk} &= - n_j n_k
  \frac{e^{-|z_j-z_k|}}{\sqrt{|z_j-z_k|}} (\widehat{z_j-z_k})_m
  \int_{\R^2} e^{x \cdot (z_k-z_j)/|z_j-z_k|} (1-a)ff'/r dx \\
  & \quad + o(e^{-R(\zb)}/\sqrt{R(\zb)}).
\end{split}
\]
Thus
\be
\la{eq:E1}
  \lan E_{\psi}^{\zb,\chi}, [\Tz_{lm}]_{\psi} \ran
  = n_l \sum_{k \not= l} c_{lk} n_k
  \frac{e^{-|z_l-z_k|}}{\sqrt{|z_l-z_k|}}
  (\widehat{z_l-z_k})_m
  + o(e^{-R(\zb)}/\sqrt{R(\zb)}).
\end{equation}
The computation of $\lan E_A^{\zb,\chi}, [\Tz_{lm}]_A \ran$ is similar,
but simpler.  We just report the result:
\be
\la{eq:E2}
  \lan E_A^{\zb,\chi}, [\Tz_{lm}]_A \ran =
  n_l \sum_{k \not= l} c'_{lk} n_k
  \frac{e^{-|z_l-z_k|}}{\sqrt{|z_l-z_k|}}
  (\widehat{z_k-z_l})_m
  + o(e^{-R(\zb)}/\sqrt{R(\zb)}).
\end{equation}
Combining~(\ref{eq:E1}) and~(\ref{eq:E2}) with~(\ref{eq:force})
yields~(\ref{eq:forceasy}).
%Due to~(\ref{eq:force}),
%the l.h.s. differs from
%$\nabla_{z_{lm}} W(\zb)$ by $o(\e)$.
%The proposition now follows from the asymptotics~(\ref{eq:decay}).
%and the relation
%$j^{(n)} = n(f^{(n)})^2\frac{1}{r}(1-a^{(n)})J \hat{x}$,
%we have
%$n\frac{1}{r}(1-a^{(n)}) = \beta_n K_1(r)(1 + o(1))$,
%and the r.h.s. is
%$-(\pi/\beta_{n_l}) \nabla_{z_{lm}} W(\zb)$
%to leading order in $\e$.
$\Box$

\begin{rem}
\la{rem:type1}
Similar computations can be made for the
Type-I case, $\lam < 1/2$.  In this case,
\[
  W(\zb) = O(e^{-m_{\lam} R(\zb)})
\]
as $R(\zb) \ra \infty$,
and the inter-vortex forces are attractive.
\end{rem}

%-----------------------------------------------------
\subsection{Approximate zero-mode property}
\la{sec:zero}

{\bf Proof of Lemma~\ref{lem:zero}.}
Set $\Lz := \E_{GL}''(\vz)$.
For any, $j$, we may write
\[
  \Lz = L_j + V_{(j)}
\]
where
$L_j := \E_{GL}''(g_{\chi_{(j)}} u^{\nj}(\cdot-z_j))$,
$\chi_{(j)} := \chi + \sum_{k \not= j} \h(\cdot-z_k)$,
and $V_{(j)}$ is
a multiplication operator satisfying
\[
  |V_{(j)}(x)| \leq c e^{-min_{k \not=j} |x-z_k|}.
\]
The notation $g_\g u$ stands for the result of acting
on $u$ by a gauge transformation $\g$.
Recall the translational modes $\Tz_{jk}$ are given
in~(\ref{eq:texp}). Using the fact that
\[
  L_j ( e^{i \chi_{(j)}} (\nabla_{A_k} \psi)_j,
  B_j \hat{e}^{\perp}_j ) = 0,
\]
we get the easy estimates
$\|L_j \Tz_{jk}\|_2 \leq ce^{-R(\zb)}$, and
\[
  \| V_{(j)} \Tz_{jk} \|_2 \leq
  c e^{-R(\zb)}.
\]
Thus
\be
\la{eq:kergl1}
  \| \Lz \Tz_{jk} \|_2 \leq c \e \log^{1/2}(1/\e).
\end{equation}

To deal with the gauge modes,
$\Gz_\g := \lan \g, \p_{\chi} \vz \ran$,
we use~(\ref{eq:gexp}), which gives
\[
  \Lz \Gz_{\g} =
  (i \g [\E_{GL}'(\vz)]_\psi, 0)
\]
and so
\be
\la{eq:kergl2}
  \| \Lz \Gz_{\g} \|_2 \leq c \e \log^{1/4}(1/\e) \| \g \|_2.
\end{equation}
Now by~(\ref{eq:Gamma1}),
and~(\ref{eq:tan1'})-~(\ref{eq:Fdef}),
$Ran Q_\s = T_{w_\s} M_{mv}$ consists of vectors of the form
$(\al \cdot \Tz + \Gz_\g, O_{L^2}(\sqrt{\e}))$
with $\al \in \R^{2m}$ and $\g \in H^1$.
This, together with
\[
  L_\s = \left( \begin{array}{cc}
  \Lz & 0 \\ 0 & {\bfone}
  \end{array} \right),
\]
and~(\ref{eq:kergl1}) and~(\ref{eq:kergl2}),
yields Lemma~\ref{lem:zero}.
$\Box$

%--------------------------------------------------
\subsection{Coercivity of the Hessian}
\la{sec:coer}

{\bf Proof of Lemma~\ref{lem:coer}}.
Suppose $\eta := (\eta_\psi, \eta_A)$
is orthogonal to each approximate
translational zero-mode, $\Tz_{jk}$, and
to the approximate gauge zero-modes,
$\Gz_\g$ (which means
$Im(\overline{\psiz} \eta_\psi) = \nabla \cdot \eta_A$
by an integration by parts).
Set $L := \Lz$.  Our first goal is to show
\[
  \lan \eta, L \eta\ran \;\; \geq \;\; c_1 \| \eta \|_{H^1}^2.
\]

Let $\{ \chi_j \}$ be a partition of unity associated to the
vortex centres.  That is, $\sum_{j=0}^m \chi_j^2 = 1$,
$\chi_j$ is supported in a ball of fixed radius
about $z_j$ ($j = 1,\ldots,m$), and $\chi_0$
is supported away from all the vortices.
By the IMS formula~(\cite{cfks}),
\[
   L = \sum \chi_j L \chi_j
   -2\sum |\nabla \chi_j|^2.
\]
We can choose $\{ \chi_j \}$ such that
$|\nabla \chi_j| \leq c R^{-1}$,
where $R := R(\zb)$.  As in Section~\ref{sec:zero}, set
\[
  L_j := \E_{GL}''(g_{\chi_{(j)}} u^{\nj}(\cdot-z_j)),
\]
and write, for each $1 \leq j \leq m$,
$L = L_j + V_{(j)}$. Since
\[
  |V_{(j)}(x)| \leq c \sum_{k \not= j} e^{-|x-z_k|},
\]
we can choose $\{ \chi_j \}$ so that
$\| V_{(j)} \chi_j \|_{\infty} \leq c \sqrt{\e}$, and so
\[
  \lan\chi_j \eta, L \chi_j \eta\ran \geq
  \lan\chi_j \eta, L_j \chi_j \eta\ran - c\sqrt{\e} \|\eta\|_2^2.
\]
for $1 \leq j \leq m$.
Also, since $\chi_0$ is supported
away from all the vortices,
\[
  \lan\chi_0 \eta, L \chi_0 \eta\ran \geq
  c_2 \|\chi_0 \eta\|_{H^1}^2
\]
for some $c_2 > 0$. Thus
\[
  \lan\eta, L \eta\ran \quad \geq \quad
  \sum_{j=1}^m \lan\chi_j \eta, L_j \chi_j \eta\ran
  + c_2 \|\chi_0 \eta\|_{H^1}^2
  - c(\sqrt{\e} + R^{-2}) \|\eta\|_{H^1}^2.
\]
Now let $\{ \tilde{T}_{jk} \}$ ($k=1,2$) be the exact translational
zero-eigenfunctions of $L_j$
(see \cite{gs} for a discussion).  We have
\[
  |\lan \tilde{T}_{jk}, \chi_j \eta \ran| \leq c\e,
\]
and
\[
  Im (e^{-i \chi_{(j)}} \bar{\psi}_j \chi_j \eta_\psi)
  - \nabla \cdot (\chi_j \eta_A) = O(R^{-1}).
\]
So by the $n$-vortex stability result of~\cite{gs}
(for $n_j = \pm1$), we have
\[
  \lan\chi_j \eta, L_j \chi_j \eta\ran
  \geq c_3 \|\chi_j \eta \|_{H^1}^2 - c\e \|\eta\|_2^2,
\]
and so
\be
\la{eq:coer}
  \lan\eta, L \eta\ran \geq
  [c_4 - c(\sqrt{\e} + R^{-2})] \| \eta \|_{H^1}^2
  \geq c_1 \| \eta \|_{H^1}^2
\end{equation}
for $\e$ sufficiently small.

For the Higgs model,
the linearized operator acts as the identity on the
momentum components:
\[
  L_\s := \cH''(w_\s) = \left( \begin{array}{cc}
  \Lz & 0 \\ 0 & \bfone \end{array} \right).
\]
Observing that $Q_\s \xi \equiv 0$ implies $\Pz \xi_1 =
O(\sqrt{\e})$ (see Eqn.~(\ref{eq:Qform}) in Appendix 2), we have
\[
  \lan \xi, L_\s \xi \ran \geq
  [\g - O(\sqrt{\e})]\|\xi_1\|_{H^1}^2 + \|\xi_2\|_2^2.
\]
This proves Lemma~\ref{lem:coer} (the upper bound is
straightforward).
$\Box$

\smallskip

\noindent
{\bf Proof of Lemma~\ref{lem:co2}}:
Set $L := \Lz$.
First observe, using $L \Pz = o(1)$, that
\[
\begin{split}
  \lan L \xi, L \xi \ran
  &=\lan L^{1/2} \xi, L L^{1/2} \xi \ran
  =\lan \Pz L^{1/2} \xi, L L^{1/2} \xi \ran \\
  &\quad + \lan \bPz L^{1/2} \xi, (\Pz + \bPz) L L^{1/2} \xi \ran
  \geq (c_1-o(1)) \|\xi\|_{H_1}^2.
\end{split}
\]
% is this true ?
Now since $\|L + \Delta\|_{H^1 \ra L^2} \leq c$,
for any $0 < \delta < 1$ we have
\[
\begin{split}
  \lan L \xi, L \xi \ran
  &= \delta \lan L \xi, L \xi \ran
  + (1-\delta) \lan L \xi, L \xi \ran \\
  &\geq \delta \lan \Delta \xi, \Delta \xi \ran
  - c \delta \|\xi\|_{H^1}^2
  + (1-\delta)(c_1-o(1)) \|\xi\|_{H^1}^2
  \geq \gt \|\xi\|_{H^2}^2
\end{split}
\]
for $\delta$ and $\e$ sufficiently small.
$\Box$

%---------------------------------------------------------
\subsection{Remainder estimates for GL functional}
\label{sec:tay}

{\bf Proof of Lemma~\ref{lem:nl2}}.
For $v = (\psi,A)$, set
\[
  R_{v}(\xi) := \E_{GL}(v + \xi) - \lan\E_{GL}'(v), \xi\ran
  - \frac{1}{2} \lan \xi, \E_{GL}'' \xi\ran.
\]
Here
\[
  \E_{GL}'(v) = \left(
  -\Delta_A \psi + \lam(|\psi|^2-1)\psi,
  -curl^2 A - Im(\bar{\psi}\nabla_A \psi) \right)
\]
and $\E''_{GL}(v)$ is the Hessian
of $\E_{GL}$ at $v$
(which we don't write out explicitly here).
After some computation, we find, for $\xi = (\eta,\al)$,
\[
\begin{split}
  R_v(\xi) &=
  \int \left\{ |\al|^2 Re(\bar{\eta} \psiz) - \al \cdot
     Im(\bar{\eta} \nabla_{\Az} \eta) \right. \\
  &\quad + \left. \lam Re(\overline{\psiz}\eta)|\eta|^2
  + \frac{1}{2} |\al|^2|\eta|^2
  + \frac{\lam}{4} |\eta|^4 \right\}
\end{split}
\]
and so using H\"older's inequality,
and the Sobolev embedding
$\| g \|_p \leq c_p \| g \|_{H^1}$
in two dimensions, we obtain easily
\be
\la{eq:tayest}
  | R_{v}(\xi) | \leq c(\|\xi\|_{H^1}^3 + \|\xi\|_{H^1}^4).
\end{equation}
$\Box$

\smallskip

\noindent
{\bf Proof of Lemma~\ref{lem:nl1}}.
The most problematic term in $N_v(\xi)$ is of the
form $\xi \nabla \xi$, so we will just
bound this one (the rest are straightforward):
\[
\begin{split}
  \| \xi \nabla \xi \|_{H^{-s}}
  &= \sup_{\| \eta \|_{H^s} = 1} |(\eta, \xi \nabla \xi)|
  \leq \sup \| \eta \xi \|_2 \| \nabla \xi \|_2 \\
  &\leq c \sup \|\eta\|_p \|\xi \|_q \| \xi \|_{H^1}
  \leq c \| \xi \|_{H^1}^2
\end{split}
\]
where $1/p + 1/q = 1/2$ and $q$ is taken large enough
so that $H^s \subset L^p$.
$\Box$

%%%%%%%%%%%%%%%%%%%%%%%%%%%%%%%%%%%%%%%%%%%%%%%%%%%%%%
\section{Appendices}

%------------------------------------------------------
\subsection{Appendix 1: Operators $V_\s$}
\la{sec:UV}

In this section we consider the key operators
$V_\s := \Lam_\s \J^{-1} \G_\s$, where $\G_\s$ and $\Lam_\s$ are
given in~(\ref{eq:Gamma2}) and~(\ref{eq:ladef'}),
and we show for them the relation~(\ref{eq:Vs}), implying, in
particular, the invertibility of $V_\s$.
We also prove the
auxiliary properties~(\ref{eq:gamap}) and~(\ref{eq:lamap}) of the
operators $\G_\s$ and $\Lam_\s$.
To this end, we use the following explicit expressions for the
basis vectors~(\ref{eq:tan}) for the tangent
space $T_{w_\s} M_{mv}$:
\be
\la{eq:tan1'}
  \ta^z_{jk} = (\Tz_{jk}, S^\s_{jk}),
\end{equation}
\be
\la{eq:tan2'}
  \ta^p_{jk} = (0, \Tz_{jk}),
\end{equation}
\be
\la{eq:tan3'}
  \ta^\chi_x = (\Gz_{\del_x}, F^\s_{\del_x}),
\end{equation}
\be
\la{eq:tan4'}
  \ta^\z_x = (0, \Gz_{\del_x})
\end{equation}
where
\be
\la{eq:Sdef}
  S^\s_{jk} :=
  ( (\p_{z_{jk}} + iA_{jk})[\phi_\s]_\pi,
  \p_{z_{jk}} [\phi_\s]_E )
\end{equation}
and
\be
\la{eq:Fdef}
  F^\s_{\del_x} := \p_{\chi(x)} \phi_\s
  = ( i \del_x [\phi_\s]_\pi, 0).
\end{equation}

In what follows we omit the super- and sub-indices $(\zb,\chi)$ and
$\s$. Using Equations~(\ref{eq:tan1'})-~(\ref{eq:Fdef}) it is not
difficult to verify properties~(\ref{eq:gamap})
and~(\ref{eq:lamap}).  For example, to show~(\ref{eq:lamap}) we
calculate using definitions~(\ref{eq:ladef'})
and~(\ref{eq:tan1'})-~(\ref{eq:Fdef}),
\[
  \Lam_\s \xi = (\al,\g,\beta,\eta),
\]
where $\xi = (\xi_1, \xi_2)$
\[
  \al_{jk} := \lan T_{jk}, \xi_1 \ran
  + \lan S_{jk}, \xi_2 \ran,
\]
\[
  \g(x) := \lan G_{\del_x}, \xi_1 \ran
  + \lan F_{\del_x}, \xi_2 \ran,
\]
\[
  \beta_{jk} := \lan T_{jk}, \xi_2 \ran,
  \quad
  \eta(x) := \lan G_{\del_x}, \xi_2 \ran.
\]
Clearly $| \al | \leq c \| \xi \|_{X_{r,s}}$
for any $r,s$, and similarly for $\beta$.
Furthermore, due to~(\ref{eq:gexp}),
\[
  \lan G_{\del_x}, \xi_1 \ran =
  Im(\bar{\psi} \xi_\psi) - div \xi_A
\]
and due to~(\ref{eq:Fdef}),
\[
  \lan F_{\del_x}, \xi_2 \ran =
  Im( \bar{\phi}_\pi \xi_\pi ).
\]
Recall that $\phi := \sum_{j=1}^m p_j \cdot \Tz_j + \Gz_\z$ and
that $\z \in H^1.$ Using that
$H^1(\R^2) \cdot H^s(\R^2) \subset H^{r'}(\R^2)$
for $r' < min(s,1),$ we obtain
\[
  \| \g \|_{H^{r-1}} \leq
  c( \|\xi_\psi\|_{H^{r-1}} + \|\xi_A\|_{H^r}
  + \|\xi_\pi\|_{H^{s}})
  \leq c \| \xi \|_{H^r \times H^{s}},
\]
provided $r-1 < min(s,1)$. Similarly, we have
\[
  \|\eta\|_{H^{s-1}} \leq
  c \| \xi_2 \|_{H^s} \leq
  c \| \xi \|_{H^{any} \times H^s}.
\]
Summing this up, we conclude that
\[
  \| \Lam_\s \xi \|_{Y_{r-1,s-1}} \leq
  c \| \xi \|_{X_{r,s}}
\]
provided $r-1 < min(s,1)$, which implies~(\ref{eq:lamap}).
(\ref{eq:gamap}) is obtained similarly.

Now we use
explicit expressions~(\ref{eq:tan1'})-~(\ref{eq:tan4'}) for the
basis~(\ref{eq:tan}) in order to establish
equation~(\ref{eq:Vs}). Equation~(\ref{eq:Vs}) follows from the
relations $\lan \ta, \J^{-1} \ta \ran = 0$,
where $\ta = \ta^z_{ij}$, $\ta^p_{ij}$, $\ta^\chi_x$, or
$\ta^\z_x$, and the relations
\be
\la{eq:skew1}
  \lan \ta^z_{ij}, \J^{-1} \ta^p_{kl} \ran
  = -\lan \ta^p_{ij}, \J^{-1} \ta^z_{kl} \ran
  = \g_{n_i} \del_{ik} \del_{jl} + O(\e\log^{1/2}(1/\e)),
\end{equation}
\be
\la{eq:skew2}
  \lan \ta^z_{ij}, \J^{-1} \ta_y^\chi \ran,
  \lan \ta^\chi_x, \J^{-1} \ta_{kl}^z \ran
 = O(\sqrt{\e}),
\end{equation}
\be
\la{eq:skew3}
\begin{split}
  \lan \ta^z_{ij}, \J^{-1} \ta_y^\z \ran &=
  \lan \ta^\z_x, \J^{-1} \ta_{kl}^z \ran =
  \lan \ta^p_{ij}, \J^{-1} \ta_y^\chi \ran \\
  &= \lan \ta^\chi_x, \J^{-1} \ta_{kl}^p \ran
  = \lan \ta^p_{ij}, \J^{-1} \ta_y^\z \ran \\
  &= \lan \ta^\z_x, \J^{-1} \ta_{kl}^p \ran = O(\e \log^{1/4}(1/\e)),
\end{split}
\end{equation}
\be
\la{eq:skew4}
  \lan \ta^\chi_x, \J^{-1} \ta_y^\z \ran
  =-\lan \ta^\z_x, \J^{-1} \ta_y^\chi \ran
  = -K_{xy}
\end{equation}
where $K_{xy}$ is the integral kernel of the operator
$K = -\Delta + |\psiz|^2$.

We will not present here proofs of all
the relations~(\ref{eq:skew1})-~(\ref{eq:skew4}),
but rather illustrate
our arguments by establishing two of the relations, say $\lan
\ta^z_{ij}, \J^{-1} \ta_y^\chi \ran$ and $\lan \ta^z_{ij}, \J^{-1}
\ta_y^\z \ran$ (see~(\ref{eq:skew2}) and~(\ref{eq:skew3})). In
what follows, we omit the superscripts in $\Tz_{jk}$,
$\Gz_{\del_x}$, $S^\s_{jk}$, and $F^\s_{jk}$, and the subscripts
in $\psiz$ and $\phi_\s$. Using Equations~(\ref{eq:tan1'})
and~(\ref{eq:tan3'}), we obtain
\[
  \lan \ta^z_{jk}, \J^{-1} \ta_y^\chi \ran
  = \lan T_{jk}, F_{\del_y} \ran
  - \lan S_{jk}, G_{\del_y} \ran.
\]
Using the explicit expressions~(\ref{eq:texp}),
(\ref{eq:gexp}), (\ref{eq:Sdef}), and~(\ref{eq:Fdef})
for the vectors on the r.h.s, we compute
\[
  \lan \ta^z_{jk}, \J^{-1} \ta^\chi_y \ran
  = -Im(e^{-i \chi} \bar{\psi}_{(jk)} \phi_\pi)
  + Im[\overline{(\p_{z_{jk}} + iA_{jk})\phi_\pi} \psi]
  + div(\p_{z_{jk}} \phi_A),
\]
where
$\psi_{(jk)}(x) = [\prod_{l \not= j} \psi^{(n_l)}(x-z_l)]
([\nabla_A]_k \psi)^{(n_j)}(x-z_j)$.
Recalling the definition~(\ref{eq:phidef}) of $\phi$,
and using $|p|+\|\z\|_{H^1} < \e$,
we conclude that~(\ref{eq:skew2}) is true.

To prove that
$\lan \ta_{ij}^z, \J^{-1} \ta_y^\z \ran = O(\e)$,
use Equations~(\ref{eq:tan1'}) and~(\ref{eq:tan3'})
to obtain
\[
  \lan \ta^z_{jk}, \J^{-1} \ta_y^\z \ran
  = \lan T_{jk}, G_{\del_y} \ran.
\]
Now using~(\ref{eq:texp}) and~(\ref{eq:gexp}) and
the equation
$curl B^{(n)} = Im(\bar{\psi}^{(n)} \nabla_{A^{(n)}} \psi^{(n)})$
and estimate~(\ref{eq:decay}), we find
\[
\begin{split}
  |\lan T_{jk}, G_{\del_y} \ran| &=
  |Im(\prod_{l \not= j}
  \overline{\psi^{(n_l)}(x-z_l)([\nabla_A]_k \psi)^{(n_j)}(x-z_j)}
  \prod_m \psi^{(n_m)}(x-z_m)) \\
  & \quad - curl B^{(n_j)}(x-z_j)| \\
  &= |Im[(\prod_{l \not= j} |\psi^{(n_l)}(x-z_l)|^2-1)
  \overline{([\nabla_A]_k \psi)^{(n_j)}(x-z_j)} \psi^{(n_j)}(x-z_j)]|
\end{split}
\]
and using Lemma~\ref{lem:*}, we see
$|\lan T_{jk}, G_{\del_y} \ran| \leq c \e \log^{1/4}(1/\e)$.

%---------------------------------------------------------
\subsection{Appendix 2: Proof of Lemma ~\ref{lem:prec}}
\la{sec:proj}

We first note that using $Q_\s \xi = 0$ and therefore
$Q_\s \p_t \xi = [Q_\s,\p_t]\xi$, and
using~(\ref{eq:Q'}) and~(\ref{eq:h'exp}), we obtain
\[
  \lan \cH'(w_\s), Q_\s \p_t \xi \ran
  = \lan \phi_\s, ([Q_\s, \p_t] \xi)_2 \ran
  + O( \e \log^{1/4}(1/\e)\|\xi\|_X \| \dot \s \|_{Y_{1,0}} )
\]
Using the above expressions for $V_\s$, $\Gamma_\s$, and
$\Lambda_\s$, and differentiating $Q_\s$ with respect to $t$ leads to
%\[
%\begin{split}
%  Q_\s \p_t \xi &= -\Omega_s^{-1}
%  (\lan \J \p_t \ta^z, \xi \ran,
%   \lan \J \p_t \ta^p, \xi \ran,
%   \lan \J \p_t \ta^\chi, \xi \ran,
%   \lan \J \p_t \ta^\z, \xi \ran) \\
%  &= -\Omega_0^{-1}
%  (\lan \J \p_t \ta^z, \xi \ran,
%   \lan \J \p_t \ta^p, \xi \ran,
%   \lan \J \p_t \ta^\chi, \xi \ran,
%   \lan \J \p_t \ta^\z, \xi \ran) \\
%  & \quad + O(\sqrt{\e} \|\xi\|_X
%  (|\dot{\zb}|+\|\p_t^{\zb} \chi\|_{H^1}
%  +|\dot{\pb}|+\|\p_t \z\|_2)).
%\end{split}
%\]
%Hence
\[
\begin{split}
  ([Q_\s, \p_t] \xi)_2 &= -
  D^{-1}_{jk,lm}
  \lan (\dot{z}_{qr}
  [\p_{z_{qr}} + \left( \begin{array}{cc}
  i A_{qr} & 0 \\ 0 & 0 \end{array} \right)]
  + \lan \p_t^{\zb} \chi, \p_\chi \ran) \Tz_{jk}, \xi_2 \ran \Tz_{lm} \\
  & \quad -
  G_{(-\Delta + |\psiz|^2)^{-1}
  \lan (\dot{z}_{qr}
  [\p_{z_{qr}} + \left( \begin{array}{cc}
  i A_{qr} & 0 \\ 0 & 0 \end{array} \right)]
  + \lan \p_t^{\zb} \chi, \p_\chi \ran) \Gz_\del, \xi_2 \ran} \\
  & \quad +
  O(\|\xi\|_X (|\dot{\pb}|+\|\p_t \z\|_2
  + \sqrt{\e}(|\dot{\zb}|+\|\p_t^{\zb} \chi\|_{H^1}))).
\end{split}
\]
Using $\phi_\s = p_{st} \Tz_{st} + \Gz_\z$, we find
\[
  \lan \phi_\s, [Q_\s \p_t \xi]_2 \ran =
  \lan \dot{z}_{qr} S^\s_{qr}
  - F^\s_{\p_t^{\zb} \chi}, \xi_2 \ran
  + O(\sqrt{\e} \|\xi\|_X (|\dot{\pb}|+\|\p_t \z\|_2
  + \sqrt{\e}(|\dot{\zb}|+\|\p_t^{\zb} \chi\|_{H^1})))
\]
and so conclude
\[
  |\lan \cH'(w_\s), Q_\s \p_t \xi \ran
  - \lan \dot{z}_{qr} S^\s_{qr}
  - F^\s_{\p_t^{\zb} \chi}, \xi_2 \ran|
  \leq c( \sqrt{\e} \|\xi\|_X (|\dot{\pb}|+\|\p_t \z\|_2
  + \sqrt{\e}\log^{1/4}(1/\e)(|\dot{\zb}|+\|\p_t^{\zb} \chi\|_{H^1}))).
\]

We turn now to computation of $\lan \cH'(w_\s), Q_\s \J L_\s \xi
\ran$. To do this, we have to refine the above computations, and
compute $Q_\s$ up to $O(\e)$. We find
\be
\la{eq:Qform}
  Q_\s = \left( \begin{array}{cc}
  \Pz & 0 \\ 0 & \Pz
  \end{array} \right) +
  \left( \begin{array}{cc}
  0 & Q_{12} \\ Q_{21} & 0
  \end{array} \right) + O(\e\log^{1/2}(1/\e))
\end{equation}
where $\Pz$ denotes the orthogonal projection onto the span of the
vectors $\Tz_{jk}$ and $\Gz_\g$, and $Q_{12}$ and $Q_{21}$ are
$O(\sqrt{\e})$. We will need the explicit form of $Q_{21}$:
\[
\begin{split}
  Q_{21} &=
  D^{-1} \lan \Tz, \cdot \ran S^\s +
  (-\Delta + |\psiz|^2)^{-1} \lan \Gz, \cdot \ran F^\s \\
  & \quad - D^{-1} \lan S^\s, \cdot \ran \Tz
  - (-\Delta + |\psiz|^2)^{-1} \lan F^\s, \cdot \ran G_\s.
\end{split}
\]
We have
\[
\begin{split}
  \lan [\cH'(w_\s)]_2, [Q_\s \J L_\s \xi]_2 \ran &=
  - \lan Q_\s \J \cH'(w_\s), L_\s \xi \ran \\
  &= \lan [Q_\s ( \phi_\s, -\E'(\vz))]_1,
  \E''(\vz) \xi_1 \ran +
  \lan [Q_\s \phi_\s, -\E'(\vz))]_2,
  \xi_2 \ran \\
  &= \lan \Pz \phi + O(\e^{3/2}\log^{1/2}(1/\e)), \E''(\vz) \xi_1 \ran \\
  &\quad + \lan -\Pz \E'(\vz) + Q_{21} \phi + O(\e^{3/2}\log^{1/2}(1/\e)),
  \xi_2 \ran \\
  &= \lan \E''(\vz) \Pz \phi, \xi_1 \ran
  - \lan \E'(\vz), \Pz \xi_2 \ran
  + \lan Q_{21} \phi, \xi_2 \ran \\
  &\quad + O(\e^{3/2}\log^{1/2}(1/\e)\|\xi\|_X).
\end{split}
\]
Now use $\E''(\vz) \Pz = O(\e\log^{1/2}(1/\e))$ and the fact that
\[
  0 = Q_\s \xi = (\Pz {\bf 1} + O(\sqrt{\e})) \xi
\]
which implies $\Pz \xi_j = O(\sqrt{\e} \|\xi\|_X)$ to find
\[
  \lan \cH'(w_\s), Q_\s \J L_\s \xi \ran
  = \lan Q_{21} \phi, \xi_2 \ran
  + O(\e^{3/2}\log^{1/2}(1/\e)\|\xi\|_X).
\]
From the form of $Q_{21}$ given above, and the fact that $\lan
\Tz, \xi_2 \ran$ and $\lan \Gz, \xi_2 \ran$ are
$O(\sqrt{\e}\|\xi\|_X)$, we see
\[
\begin{split}
  \lan \cH'(w_\s), Q_\s \J L_\s \xi \ran
  &= \lan D^{-1} \lan \Tz, \phi_\s \ran S^\s
  + (-\Delta + |\psiz|^2)^{-1}
  \lan \Gz, \phi_\s \ran F^\s, \xi_2 \ran \\
  &\quad + O(\e^{3/2}\log^{1/2}(1/\e)\|\xi\|_X) \\
  &= \lan p_{jk} S^\s_{jk} - F_\z, \xi_2 \ran
  + O(\e^{3/2}\log^{1/2}(1/\e)\|\xi\|_X).
\end{split}
\]
Combining this with the above computation of $\lan \cH'(w_\s),
Q_\s \p_t \xi \ran$, with the facts that $|\dot \pb| = |\dot \pb +
\nabla_z W(\zb)| + O(\e)$ and $|\dot \zb| + \| \p_t^{\zb} \chi
\|_{H^1} = |\dot \zb - \pb| + \| \p_t^z \chi - \z\|_{H^1} +
O(\sqrt{\e})$ proves Lemma~\ref{lem:prec} $\Box$

%---------------------------------------------------
\subsection{Appendix 3: Two technical lemmas}
\la{sec:app3}

\begin{lem}
\la{lem:*} Let $0 < \al \leq \beta$ and $0 \leq \del, \gamma <
3/2$. Then \be \la{eq:lem1}
  \int_{\R^2} \frac{e^{-\al|x|}e^{-\beta|x-a|}}
  {|x|^{\gamma}|x-a|^{\del}} dx
  \leq c \frac{e^{-\al|a|}}{|a|^{\gamma+\del-2}}
  \left\{ \begin{array}{cc}
  |a|^{-1/2} & \al=\beta \\
  |a|^{\del-2} & \al < \beta
  \end{array} \right. .
\end{equation}
\end{lem}

\noindent
{\bf Proof:}
We prove only the case $\alpha=\beta$, since
the remaining cases follow from Lemma~\ref{lem:**} below.
Define $I$ by
\[
\begin{split}
  2 |a|^{2-\gamma-\del} I &:=
  \int_{\R^2} \frac{e^{-\al|x|-\beta|x-a|}}
  {|x|^{\gamma}|x-a|^{\del}} dx \\
  &= \int_0^\infty \frac{dr}{r^{\gamma-1}}
  \int_0^{2 \pi} d \theta
  \frac{e^{-\al r - \beta \sqrt{r^2+|a|^2-2r|a|\cos(\theta)}}}
  {(r^2+|a|^2-2r|a|\cos(\theta))^{\del/2}}.
\end{split}
\]
Changing variables to $u = 1-\cos(\theta)$ and $t=r/|a|$,
and using $\alpha=\beta$, we estimate
\be
\la{eq:lem2}
\begin{split}
  I &\leq \int_0^{\infty} \frac{dt}{t^{\gamma-1}}
  \int_0^1 \frac{du}{\sqrt{u}}
  \frac{e^{-\alpha|a|(t + \sqrt{(t-1)^2+2tu})}}
  {((t-1)^2+2tu)^{\delta/2}} \\
  &= \int_0^{1/2} \int_0^1 + \int_{1/2}^1 \int_0^1
  + \int_1^\infty \int_0^1
  =: I_1 + I_2 + I_3.
\end{split}
\end{equation}
Using $(t-1)^2+2tu \geq [1-t+tu]^2$ (for $0 \leq u \leq 1$),
we have
\be
\la{eq:lem3}
\begin{split}
  I_1 &\leq c \int_0^{1/2} \frac{dt}{t^{\gamma-1}}
  \int_0^1 \frac{du}{\sqrt{u}} e^{-\alpha|a|(1+tu)} \\
  &= c e^{-\alpha|a|}/\sqrt{|a|}
  \int_0^{1/2} \frac{dt}{t^{\gamma-1/2}}
  \int_0^{|a|t} e^{-\alpha v}/\sqrt{v} dv \\
  &\leq c e^{-\alpha |a|}/\sqrt{|a|}.
\end{split}
\end{equation}
Now estimating $\sqrt{(t-1)^2+2tu} \ge \sqrt{(1-t)^2+u}$ for $t
\ge 1/2$ and changing the variables of integration as $x=|a|(1-t)$
and $y=|a|\sqrt{u}$ we obtain
\[
\begin{split}
  I_2 &\leq c e^{-\alpha|a|}/|a|^{2-\del}
  \int_0^{|a|} dx \int_0^{|a|} dy
  \frac{e^{-\alpha(-x+\sqrt{x^2+y^2})}}{(x^2+y^2)^{\del/2}} \\
  &\leq c e^{-\alpha|a|}/|a|^{2-\del}
  \int_0^{|a|} dx \int_0^{|a|} dy
  \frac{e^{-\alpha y^2/(2\sqrt{x^2+y^2})}}{(x^2+y^2)^{\del/2}} \\
  &\leq e^{-\alpha|a|}/|a|^{2-\del}
  \int_0^{2|a|} \frac{dr}{r^{\del-1}}
  \int_0^{2\pi} d \theta e^{-\alpha r \cos^2(\theta)/2}.
\end{split}
\]
We have
\[
\begin{split}
  \int_0^{2\pi} d \theta e^{-\alpha r \cos^2(\theta)/2}
  &= c \int_{-1}^1 \frac{ds}{\sqrt{1-s^2}}
  e^{-r(1+s)\alpha/4} \\
  &\leq c \int_0^1 \frac{ds}{\sqrt{1-s}}e^{-r(1-s)\alpha/4}
  + c e^{-r \alpha/4} \\
  &= \frac{c}{\sqrt{r}} \int_0^r \frac{dv}{\sqrt{v}}
  e^{-\sqrt{\alpha}/4} + c e^{-r \alpha/4}
  \leq c/\sqrt{r}.
\end{split}
\]
Hence
\be
\la{eq:lem4}
  I_2 \leq c \frac{e^{-\alpha|a|}}{\sqrt{|a|}}.
\end{equation}
Finally,
\begin{equation}
\la{eq:lem5}
\begin{split}
  I_3 &\leq \int_0^\infty \frac{d \tau}{(\tau+1)^{\gamma-1}}
  \int_0^1 \frac{du}{\sqrt{u}}
  \frac{e^{-\alpha|a|(\tau+1+\sqrt{\tau^2+2u})}}
  {(\tau^2+2u)^{\del/2}} \\
  &= 2 e^{-\alpha|a|} \int_0^\infty
  \frac{d \tau}{(\tau+1)^{\gamma-1}}
  \int_0^1 dv \frac{e^{-\alpha|a|(\tau+\sqrt{\tau^2+2v^2})}}
  {(\tau^2+2v^2)^{\del/2}} \\
  &\leq c e^{-\alpha|a|}/|a|^{2-\delta}.
\end{split}
\end{equation}
Estimates~(\ref{eq:lem2})-~(\ref{eq:lem5})
imply~(\ref{eq:lem1}).
$\Box$

\begin{lem}
\la{lem:**} Suppose $b(x)$ is a function satisfying $|b(x)| \leq c
e^{-m|x|}$ for some $m>1$, and $e(x)$ is a bounded function with
asymptotic behaviour $e(x) = c_1 e^{-|x|}/\sqrt{|x|}(1 +
O(1/|x|))$ as $|x| \ra \infty$.  Fix $z \in \R^2$ and set
\[
  I(z) := \int_{\R^2} b(x) e(x-z) dx.
\]
Then
\be
\la{eq:asy}
  I(z) = c_1 e^{-|z|}/\sqrt{|z|}
  \int_{\R^2} e^{x \cdot z/|z|} b(x) dx [1 + O(1/|z|)]
\end{equation}
as $|z| \ra \infty$.
\end{lem}

\noindent
{\bf Proof}:
Choose $\alpha$ with $1/m < \alpha < 1$.
Let $D_{\alpha |z|}$ denote the disk of radius
$\alpha |z|$ about the origin.  We have
\be
\la{eq:out}
  |\int_{\R^2 \backslash D_{\alpha |z|}} b(x) e(x-z) dx|
  \leq c \int_{\alpha |z|}^\infty e^{-mr} r dr
  \leq c |z| e^{-m \alpha |z|} << e^{-|z|}/|z|^p
\end{equation}
for any $p$.
On $D_{\alpha |z|}$, $|x| << |z|$, and we Taylor expand:
\[
  |z-x|^{-1/2} = |z|^{-1/2}[1 + O(|x|/|z|)],
\]
\[
  e^{-|z-x|} = e^{-|z|}e^{x \cdot z/|z|}[1 + O(|x|^2/|z|)],
\]
yields
\be
\la{eq:in}
  \int_{D_{\alpha |z|}} b(x) e(x-z) dx
  = e^{-|z|}/\sqrt{|z|} \int_{D_{\alpha |z|}}
  e^{x \cdot z/|z|} b(x) [1+O(|x|^2/|z|)].
\end{equation}
Estimates~(\ref{eq:out}) and~(\ref{eq:in})
together yield~(\ref{eq:asy}).
$\Box$

%%%%%%%%%%%%%%%%%%%%%%%%%%%%%%%%%%%%%%%%%%%%%%%%%
%%%               References                %%%%%
%%%%%%%%%%%%%%%%%%%%%%%%%%%%%%%%%%%%%%%%%%%%%%%%%

\noindent{Stephen Gustafson}, gustaf@math.ubc.ca \\
University of British Columbia,
Vancouver, BC, Canada V6T 1Z2

\bigskip

\noindent{I.M. Sigal}, \\
University of Notre Dame,
Notre Dame, IN 46556-4618, and \\
University of Toronto,
Toronto, ON, Canada M5S 3G3.

\end{document}